\documentclass[twoside,10pt]{article}

\usepackage{ifthen,ifpdf}
\ifpdf
\usepackage[pdftex]{hyperref}	\hypersetup{colorlinks=true,linkcolor=black,citecolor=black,urlcolor=blue,breaklinks=true,plainpages=false}
\usepackage[pdftex]{graphicx}
\usepackage[usenames,pdftex]{color}
\else
\usepackage[colorlinks=true,breaklinks=true,plainpages=false]{hyperref}
\usepackage{graphicx}
\usepackage[usenames]{color}
\fi

\usepackage{jmlr2e}
\usepackage{makecell}
\usepackage{amscd,amsxtra,amsfonts,amsmath,amssymb,multirow}
\usepackage{wrapfig}
\usepackage[footnotesize]{caption}
\usepackage[textwidth=0.8in,textsize=footnotesize]{todonotes}
\usepackage{pdflscape}
\usepackage{footmisc}

\newtheorem{thm}{Theorem}
\newtheorem{defn}{Definition}
\newtheorem{lem}{Lemma}[thm]

\ShortHeadings{Persistent-Homology-based Machine Learning and its Applications -- A Survey}{Pun, Xia, and Lee}
\firstpageno{1}

\begin{document}

\title{Persistent-Homology-based Machine Learning and its Applications -- A Survey}

\author{\name Chi Seng Pun\thanks{~The first two authors contribute equally.} \email cspun@ntu.edu.sg \\
		\addr School of Physical and Mathematical Sciences \\
		Nanyang Technological Technology \\
		Singapore 637371
		\AND
		\name Kelin Xia\footnotemark[1] \email xiakelin@ntu.edu.sg \\
		\addr School of Physical and Mathematical Sciences \\
		and School of Biological Sciences \\
		Nanyang Technological Technology \\
		Singapore 637371
       \AND
  	   \name Si Xian Lee \email lees0159@e.ntu.edu.sg \\
  	   \addr School of Physical and Mathematical Sciences \\
  	   Nanyang Technological Technology \\
  	   Singapore 637371  }

\editor{TBC}

\maketitle

\begin{abstract}
A suitable feature representation that can both preserve the data intrinsic information and reduce data complexity and dimensionality is key to the performance of machine learning models. Deeply rooted in algebraic topology, persistent homology (PH) provides a delicate balance between data simplification and intrinsic structure characterization, and has been applied to various areas successfully. 
However, the combination of PH and machine learning has been hindered greatly by three challenges, namely topological representation of data, PH-based distance measurements or metrics, and PH-based feature representation. With the development of topological data analysis, progresses have been made on all these three problems, but widely scattered in different literatures. In this paper, we provide a systematical review of PH and PH-based supervised and unsupervised models from a computational perspective. Our emphasizes are the recent development of mathematical models and tools, including PH softwares and PH-based functions, feature representations, kernels, and similarity models. Essentially, this paper can work as a roadmap for the practical application of PH-based machine learning tools. 
Further, we consider different topological feature representations in different machine learning models, and investigate their impacts on the protein secondary structure classification.
\end{abstract}

\begin{keywords}
  Persistent homology, machine learning, persistent diagram, persistent barcode, kernel, feature extraction
\end{keywords}

\section{Introduction}
Machine learning models have achieved tremendous success in computer vision, speech recognition, image analysis, text analysis, etc. However, their application in high-dimensional complicated systems have been hindered significantly by proper feature representations. Mathematically, features from geometric analysis can characterize the local structure information very well, but tend to be inundated with details and will result in data complexity. Features generated from traditional topological models, on the other hand, preserve the global intrinsic structure information, but they tend to reduce too much structure information and are rarely used in quantitative characterization. Recently, persistent homology (PH) has been developed as a new multiscale representation of topological features \citep{Edelsbrunner:2002,Zomorodian:2005,Zomorodian:2008}. It is able to bridge the gap between geometry and topology, and open up new opportunities for researchers from mathematics, computer sciences, computational biology, biomathematics, engineering, etc. The essential idea of PH is to employ a filtration procedure, so that each topological generator is equipped with a geometric measurement. In this filtration process, a series of nested simplicial complexes encoded with structural topological information from different scales are produced. The intrinsic topological properties can be evaluated. It is found that some topological invariants ``live" longer in these simplicial complexes, whereas others ``die" quickly when filtration value changes. The ``lifespans" or ``persisting times" of these invariants are directly related to geometric properties \citep{Dey:2008,Dey:2013,Mischaikow:2013}. Simply speaking, long-lived Betti numbers usually represent large-sized features. In this way, topological invariants can be quantified by their persistence in the filtration process. In other words, their lifespans or persistenting times give a relative geometric measurement of the associated topological properties. PH provides a very promising way of structure representation, and has been applied to various fields, including shape recognition \citep{DiFabio:2011}, network structure \citep{Silva:2005,LeeH:2012,Horak:2009}, image analysis \citep{Carlsson:2008,Pachauri:2011,Singh:2008,Bendich:2010,Frosini:2013}, data analysis \citep{Carlsson:2009,Niyogi:2011,BeiWang:2011,Rieck:2012,XuLiu:2012}, chaotic dynamics verification \citep{Mischaikow:1999}, computer vision \citep{Singh:2008}, computational biology \citep{Kasson:2007,YaoY:2009, Gameiro:2013,KLXia:2014c, KLXia:2015a,BaoWang:2016a, KLXia:2015c,KLXia:2015b}, amorphous material structures \citep{hiraoka:2016hierarchical,saadatfar:2017pore}, etc. Various softwares, including JavaPlex \citep{javaPlex}, Perseus  \citep{Perseus}, Dipha \citep{Dipha}, Dionysus \citep{Dionysus}, jHoles \citep{Binchi:2014jholes}, GUDHI \citep{gudhi:FilteredComplexes}, Ripser \citep{bauer2017ripser}, PHAT \citep{bauer2014phat}, DIPHA \citep{bauer2014distributed}, R-TDA package \citep{fasy:2014introduction}, etc, have been developed. The results from PH can be visualized by many methods, including persistent diagram (PD) \citep{Mischaikow:2013}, persistent barcode (PB) \citep{Ghrist:2008barcodes}, persistent landscape \citep{Bubenik:2007,bubenik:2015}, persistent image \citep{adams2017persistence}, etc. When structures are symmetric or have unique topological properties, topological fingerprints can be obtained from the PH analysis and further used in the quantitative characterization of their structures and functions \citep{KLXia:2014c,KLXia:2015a,xia2018persistent}. However, when the systems get more complicated, it becomes more challenging to directly build models on the PB or PD. Instead, machine learning models can be employed to extract the important information or to learn from these topological features.

The persistent-homology-based machine learning (PHML) models have been used in various areas, including image analysis \citep{bae2017beyond,han2016deep,qaiser2018fast,makarenko2016texture,obayashi2018persistence,giansiracusa2017persistent}, shape analysis \citep{bonis2016persistence,zhou2017exploring,li2014persistence,guo2018sparse,zeppelzauer2018study}, time-series data analysis \citep{seversky2016time,anirudh2016autism,zhang2015early,wang2014persistence,umeda2017time}, computational biology \citep{pachauri2011topology,dey2018protein,ZXCang:2015}, noise data\citep{niyogi2011topological}, sphere packing \citep{robins2016principal}, language analysis \citep{zhu2013persistent}, etc. Especially, topological-based machine learning has achieved remarkable success in drug design \citep{cang:2017topologynet,cang:2017integration,nguyen:2017rigidity,cang:2017analysis,cang:2018representability,wu:2018quantitative}. They have delivered the best results in 10 out of the 26 competitive tasks of D3R Grand Challenges \citep{nguyen2018mathematical}, which is widely regarded as the most difficult challenge in drug design. The huge success shows that, with the proper topological simplification, specially-designed PH models can not only preserve the critical chemical and biological information, but also enables an efficient topological description of biomolecular interactions \citep{cang:2018representability}. However, unlike geometric and traditional topological models, several problems exist in PHML models, including the construction and generation of meaningful metrics, kernels, and feature vectors.

The unique format of PH outcomes poses a great challenge for a meaningful metrics. To solve this problem, distance measurements or metrics \citep{mileyko2011probability,turner2014frechet,chazal2017robust,carriere2018metric,anirudh2016riemannian} has been considered, including Gromov-Hausdorff distance, Wasserstein distance \citep{mileyko2011probability}, bottleneck distance \citep{mileyko2011probability}, probability-measure-based distance \citep{chazal2011geometric,marchese2017signal,chazal2017robust}, Fisher information metric \citep{anirudh2016riemannian}. These metric definitions are usually based on PD, which can be considered as a two dimensional point distribution. However, unlike common point cloud data, PD points have different topological significance. A PD point closed to diagonal line, meaning its death time is only slightly larger than its birth time, is usually regarded as less ``useful" than the ones that are far away from the diagonal line. Because short persisting times indicate ``insignificant" topological features or ``noise". Although this is not true in various physical, chemical and biological data, in which short-lived invariants have important physical meanings (such as bond lengths, benzene ring size, etc), intrinsic topological features are associated the long-lived topological generators. Given these considerations, pseudo points are introduced and placed on the diagonal lines to guarantee the same number of points in two PDs for comparison. Further, the Wasserstein distance is used to measure the best possible matching from one PD to the other, and bottleneck distance is just a special type of Wasserstein distance. It is found that under Wasserstein metric, PDs is complete and separable \citep{mileyko2011probability}. Fr\'{e}chet means and and variances, which are generalization of mean and variance to a general metric space, can also be defined on PDs \citep{turner2014frechet} and further enable a probabilistic analysis of PDs.  In Wasserstein distance, cardinality difference between two PDs is not explicitly penalized, instead the unmatched points from the two PDs are penalized by their distance to the PD diagonal lines. A probability-measure-based distance function is proposed \citep{chazal2011geometric,marchese2017signal,chazal2017robust} to account for changes in small persistence and cardinality. It also enables statistical structure in the form of Fr\'{e}chet means and variances and establishes a classification scheme. Once PDs are transformed into probability density functions, Fisher information metric (FIM) can also be used to evaluate the distance between PDs \citep{anirudh2016riemannian}. All these distance measurements and metrics can be used in kernel constructions, and further used in machine learning models.

With the significant role of kernels in machine learning models, various PH-based kernels have been proposed, including persistence scale space (PSS) kernel \citep{reininghaus2015stable}, universal persistence scale space (u-PSS) kernel \citep{kwitt2015statistical}
persistence weighted Gaussian (PWG) kernel \citep{kusano2016persistence}, sliced Wasserstein (SW) kernel \citep{carriere2017sliced}, persistence Fisher (PF) kernel \citep{le2018riemannian}, geodesic topological kernels \citep{padellini2017supervised}, etc \citep{chazal2017robust,carriere2018metric}. Constructed from PDs/PBs, these kernels can be directly incorporated into kernel-based machine learning methods, including kernel perception, kernel support vector machine, principal components analysis, spectral clustering, Gaussian process, etc. In PSS kernel, points in a PD together with its counterparts, which are generated by flipping the points in the PD along the diagonal line, are associated or convolution with Gaussian kernels. A tunable scale parameter is incorporated into these Gaussian kernels, thus for each PD, a multiscale function is generated. The PSS kernel is defined as the inner product of the multiscale functions. Further, the u-PSS kernel is a generalization of PSS, so that it satisfies the universal property. The PWG kernel is designed to attenuate the influence from the short-lived topological generators, which are generally considered as noise. Computationally, random Fourier features can be employed to increase the efficiency of calculating PWG kernel. In SW kernel, points in PDs are projected on to lines passing through the origin. By systematically varying the slops of these lines, the distance between two PDs is transformed into the measurement of distances between their projected points on these lines. In PF kernel, each PD is transformed into a normalized probability density function, of which the square-root forms a unit Hilbert sphere. Based on the Hilbert sphere space, the Fisher information metric (FIM) can be used to define an inner product metric and a kernel is generated by incorporating the FIM into an exponential function. In a geodesic topological kernel, Wasserstein distance is directly incorporated into an exponential function. However, the geodesic topological kernels are not positive definite. PH-based kernels are one way to combine topological information with machine learning models, another important way is to generate unique vectors made of topological features and use them in machine learning models.

Topological features can be extracted from PDs/PBs. The simplest way of PD/PB-based feature generation is to collect their statistical properties, such as the sum, mean, variance, maximum, minimum, etc \citep{ZXCang:2015}. Special topological properties, such as the total Betti number in a certain filtration value, can also be considered as features \citep{ZXCang:2015}. Further, a set of coordinates, which reflect properties of the corresponding ring of functions, are proposed as the topological features in digit number classification \citep{adcock2016ring}. Tropical coordinates on the space of barcodes, which are stable with respect to Wasserstein distances, are also important topological features \citep{kalivsnik2018tropical}. However, all these methods use only partial information from PB.  A more systematical way of constructing topological feature vectors from PDs/PBs is the binning approach \citep{cang:2017topologynet,cang:2018representability}. The essential idea is to discrete the PB or PB into various elements, which are then concatenated into a feature vector. More specifically, for PB, the filtration time can be subdivided into equally-sized bins. The persistent Betti numbers in these bins form a vector. Similarly, for PDs, the binning of the birth and death time will subdived them into a equally-sized mesh. The number of topological generators in each grid box is an element of the topological feature vector. Mathematically, binning is just to employ a Cartesian grid discretization, which is simplest way to discretize a computational domain. Recently,  a new way of feature representation, known as persistent codebooks, has been proposed \citep{zielinski2018persistence,bonis2016persistence}. In these models, a fixed-sized vector representation of PDs is achieved by clustering and bag of words quantization methods, including bag of words encoding, vector of locally aggregated descriptor, and Fisher vector.
Further, instead of generating features directly from PDs/PBs, another type of feature construction is to construct persistent functions, discretize them into elements and concatenate them into a vector \citep{bubenik2017persistence,adams2017persistence}. These persistent functions includes persistent Betti number, persistent Betti function \citep{Xia:2018multiscale}, persistent landscape \citep{bubenik2017persistence}, persistent surface \citep{adams2017persistence}, etc. Similar to the binning approach, these functions are discretized first and then concatenated into a feature vector. Moreover, recently an image representation of PDs/PBs has been found extremely useful in drug design \citep{cang:2017topologynet,cang:2018representability}. In this model, biomolecular structure has been decomposed into element-specific models, each of which has its feature vector from binning approach. Instead of using these feature vectors directly, an image representation is proposed by stacking all these feature vectors together systematically. This image representation is essentially a multidimensional persistence and can be generalized to three dimensional volumetric data. More importantly, the input of topological features as images suits perfectly with the convolutional neural network (CNN) model. In this way, it really brings out the great power of persistence representation. Finally, a persistent path and signature feature model has been proposed \citep{chevyrev2018persistence}. In this model, the barcode information is embedded into a persistent path. And the persistent path is further transformed into tensor series, which can be represented into a feature vector.
It is worth mentioning that for PHML models, a great promise comes from new ways of topological representations that can incorporate more structure information, including persistent local homology \citep{bendich2015multi,fasy2016exploring,bendich2007inferring,bendich2012local,ahmed2014local}, element specific PH \citep{cang:2017integration,cang:2017analysis,cang:2017topologynet,cang:2018representability,wu:2018quantitative}, weighted PH \citep{ren2017weighted,wu2018weighted}, multidimensional PH \citep{Carlsson:2009theory,Carlsson:2009computing,Cohen:2006,Cerri:2013,KLXia:2015c,KLXia:2015d}, etc. From the statistical perspective, when a large-scale feature vector is used to represent the topological features from PBs/PDs, the PHML models will also suffer from the curse of dimensionality. In this situation, techniques for variable selection and regularization should be considered, such as the high-dimensional statistical tools used in different fields: \citep{fan2008,Cai2011,Pun2018,Pun2016,Chiu2017,Pun2019,HP2018} and the dropout approach in neural networks \citep{srivastava2014dropout}, etc.


In this paper, we review the persistent-homology-based machine learning (PHML) models and discuss its application in protein structure classification. Our focus is to provide a general guidance, as demonstrated in Figure \ref{fig:Flowchart}, for the practical application of topology-aware machine learning models. We also deliver a systematical study of PH-based metrics, kernels, and feature vectors. The comparison and analysis of these methods will help to generate more robust and efficient algorithms and models. In application part, we discuss the combination of different topological feature selections with machine learning models. 

The remainder of this paper is organized as follows. Section \ref{sec:general_pipeline} is a general pipeline for the PHML models. A more detailed discussion of the fundamentals of PH is presented in Section \ref{sec:PH_Data_Modeling}. Section \ref{Sec:PH_SL} is devoted for PH-based functions, metrics, kernels, and feature vectors, and the construction of PH-based unsupervised learning and supervised learning. We test the feature selection and its relation with different machine learning models by a protein classification example in Section \ref{sec:numericalstudy}. The paper ends with a conclusion.

\section{General Pipeline for Persistent-Homology-based Machine Learning}\label{sec:general_pipeline}

\begin{figure}[!ht]
	\begin{center}
		\begin{tabular}{c}
			\includegraphics[width=0.9\textwidth]{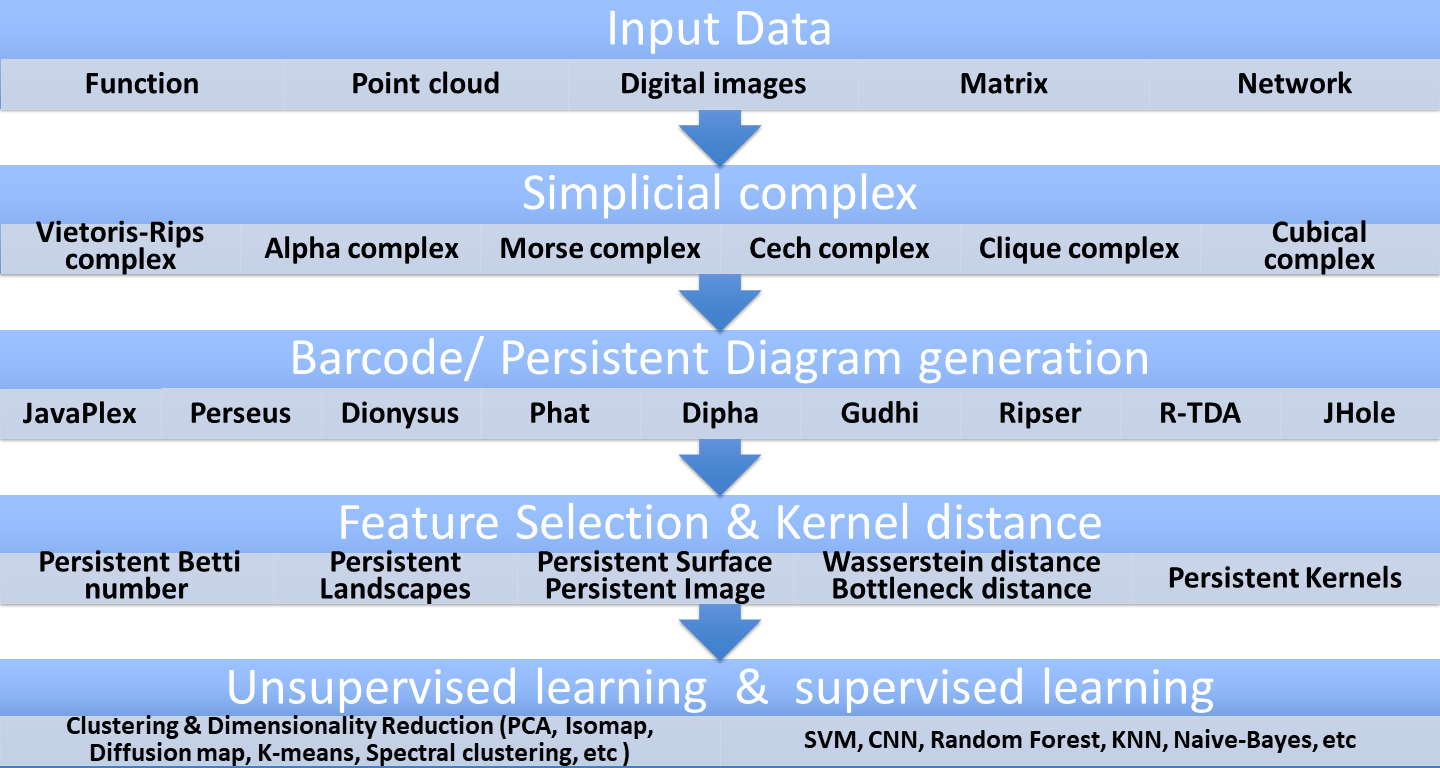}
		\end{tabular}
	\end{center}
	\caption{The flowchart for PHML modeling. For different data types,  different simplicial complexes are considered. With a suitable filtration parameter, PH analysis can be implemented by designed softwares. The result from the PH is transformed into feature vectors, distance measurements or special kernels, which are further combined with supervised or unsupervised methods from machine learning.}
	\label{fig:Flowchart}
\end{figure}

The essential idea for PHML is to extract topological features from the data using PH, and then combine these features with machine learning methods, including both supervised learning and unsupervised learning approaches. As illustrated in Figure \ref{fig:Flowchart}, PHML in computation can be divided into four steps, i.e., simplicial complex construction, PH analysis, topological feature extraction and topology-based machine learning.

The first step is to construct simplicial complex (or complexes) from the studied data. In topological modeling, we may have various types of data, including functional data, point cloud data, matrixes, networks, images, etc.  These data need to be represented by suitable simplicial complexes. Roughly speaking, the simplicial complex can be viewed as a set of combinatorial elements generated from the discretization of spaces. Cartesian grids and triangular meshes can be viewed as two dimensional simplicial complexes. Depending on the type of studied data, different simplex models can be considered, including Vietories-Rips complex, $\check{C}ech$ complex, Alpha complex, Morse complex, cubical complex, clique complex, CW complex, etc. The detailed explanation of simplicial complexes will be given in Section \ref{sec:simplicial_complex}.

The second step is the PH analysis. PH is a newly-invented model deeply-rooted in algebraic topology, computational topology and combinatorial topology. In PH \citep{Edelsbrunner:2002,Zomorodian:2005,Carlsson:2009,Edelsbrunner:2010,zomorodian2005topology}, algebraic tools, such as quotient group, homology, exact sequence, etc, are used to characterize topological invariants, including connected components, circles, rings, channels, cavities, voids, etc. Unlike traditional topological models, which capture only the intrinsic structure information and tend to ignore all geometric details, PH works as a multiscale topological representation and enables to embed geometric information back into topological invariants. This is achieved through a new idea called filtration (see Figure \ref{fig:filtration} for details). With a suitable filtration process, the persistence of topological invariants can be calculated. Roughly speaking, the persistence tells you the geometric size of the invariant. Different softwares are currently available for PH analysis of different data structures (see Figure \ref{fig:data_SimCom} for details).

The third step is to extract meaningful topological features from PH results. The PH results are usually represented as PBs or PDs (see Figure \ref{fig:barcodes_binning} for details). Neither of these representations can be used directly as input for machine learning models. Therefore, we need to transform the PH results into representations, which can be easily incorporated into machine learning models. Currently, two basic approaches are available. The first one is to generate special PH-based kernels or similarity measurements. These kernels and similarity matrixes can then be combined with PCA, SVM, K-means, spectral clustering, isomap, diffusion map, and other learning models to perform the data classification and clustering. The second one is to generate topological feature vectors from PHA results. Binning approach is the most commonly used approach to discretize the PB or PD into a feature vector. The detailed discussion will be given in Section \ref{Sec:PH_SL}.

The last step is to combine the topological features with machine learning algorithms. Essentially, these features can be used directly used in both unsupervised learning and supervised learning models. Depending on the learning models, we should consider different feature selection and representation models to bring out the best performance of the model.


\section{Persistent Homology for Data Modeling and Analysis}\label{sec:PH_Data_Modeling}

This section is devoted to the fundamentals of PH. We present a brief review of simplicial complexes and PH analysis in the first and second parts, respectively. For interested readers, a more detailed description of simplicial complex and PH can be found in papers \citep{Zomorodian:2005,Edelsbrunner:2002,Kaczynski:2004,Carlsson:2009,Mischaikow:2013,Edelsbrunner:2010}.

\subsection{Simplicial Complexes} \label{sec:simplicial_complex}

A simplicial complex is a combination of simplexes under certain rules. It can be viewed as a generalization of network or graph model.

\subsubsection{Abstract Simplicial Complex}

A simplex is the building block for simplicial complex. It can be viewed as a generalization of a triangle or tetrahedron to their higher dimensional counterparts.
\begin{defn}  A geometric $k$-simplex $\sigma^k=\{v_0,v_1,v_2,\cdots,v_k\}$ is the convex hull formed by $k+1$ affinely independent points $v_0,v_1,v_2,\cdots,v_k$ in Euclidean space $R^d$ as follows,
	$$\label{eq:couple_matrix1}
	\sigma^k=\left\{\lambda_0 v_0+\lambda_1 v_1+ \cdots +\lambda_k v_k~\left|~\sum^{k}_{i=0}\lambda_i=1;~0\leq \lambda_i \leq 1,~i=0,1, \cdots,k \right.\right\}.
	$$
	A face $\tau$ of $k$-simplex $\sigma^k$ is the convex hull of a non-empty subset. We denote it as $\tau \leq \sigma^k$
\end{defn}
Geometrically, a 0-simplex is a vertex, a 1-simplex is an edge, a 2-simplex is a triangle, and a 3-simplex represents a tetrahedron. An oriented $k$-simplex $[\sigma^k]$ is a simplex together with an orientation, i.e., ordering of its vertex set.
Simplices are the building block for (geometric) simplicial complex.

\begin{defn}\label{def:two_conditions}
	A geometric simplicial complex $K$ is a finite set of geometric simplices that satisfy two essential conditions:
	\begin{enumerate}
		\item Any face of a simplex from  $K$  is also in  $K$.
		\item The intersection of any two simplices in  $K$ is either empty or shares faces.
	\end{enumerate}
\end{defn}
The dimension of $K$ is the maximal dimension of its simplexes. A geometric simplicial complex $K$ is combinatorial set, not a topological space. However, all the points of $R^d$ that lie in the simplex of $K$ aggregate together to topologize it into a subspace of $R^d$, known as polyhedron of $K$.

Graphs and networks, which are comprised of only vertices and edges, can be viewed as a simplicial complex with only 0-simplex and 1-simplex.

\begin{defn}
	An abstract simplicial complex $K$ is a finite set of elements $\{v_0,v_1,v_2,\cdots,v_n\}$ called abstract vertices, together with a collection of subsets $(v_{i_0},v_{i_1},\cdots,v_{i_m})$ called abstract simplexes, with the property that any subset of a simplex is still a simplex.
\end{defn}

For an abstract simplicial complex $K$, there exists a geometric simplicial complex $K'$ whose vertices are in one-to-one correspondence with the vertices of $K$ and a subset of vertices being a simplex of $K'$ if and only if they correspond to the vertices of some simplex of $K$. The geometric simplicial complex $K'$ is called the geometric realization of $K$.

Based on different data types, there are various ways of defining simplicial complexes. The most commonly used ones are $\check{C}$ech complex, Vietoris-Rips complex, Alpha complex, Clique complex, Cubic complex, Morse complex, etc.

\subsubsection{$\check{C}ech$ Complex and  Vietoris-Rips Complex}
Let $X$ be a point set in Eucledian space $R^d$ and $\mathcal{U}$ is a good cover of $X$, i.e., $X\subseteq \bigcup_{i\in I}U_i$.

\begin{defn} The nerve $\mathcal{N}$ of $\mathcal{U}$ is defined by the following two conditions:
	\begin{enumerate}
		\item $\emptyset \in \mathcal{N}$;
		\item If $\bigcap_{j\in J}U_j \neq \emptyset$ for $J\subseteq I$, then $J \in \mathcal{N}$.
	\end{enumerate}
\end{defn}

\begin{thm}(Nerve theorem)
	The geometric realization of the nerve of $\mathcal{U}$ is homotopy equivalent to the union of sets in $\mathcal{U}$.
\end{thm}

We can define $B(X,\varepsilon)$ to be the closed balls of radius $\varepsilon$ centred at $x \in X$, then the union of these balls is a cover of space $X$ and $\check{C}ech$ complex is the nerve of this cover.

\begin{defn} The $\check{C}ech$ complex with parameter $\varepsilon$ of $X$ is the nerve of the collection of balls $B(X,\varepsilon)$. The $\check{C}$ech complex $\check{C}_{\varepsilon}(X)$ can be represented as $\check{C}_{\varepsilon}(X):=\left\{ \sigma \in X | \bigcap_{x \in \sigma} B(X,\varepsilon)\neq\emptyset \right\}$. 
\end{defn}

\begin{defn}The Vietoris-Rips complex (Rips complex) with parameter $\varepsilon$ denoted by $\mathcal{R}_{\varepsilon}(X)$, is the set of all $\sigma \subseteq X$, such that the largest Euclidean distance between any of its points is at most $2 \varepsilon$.
\end{defn}

It should be noticed that both $\check{C}ech$ complex and Vietoris-Rips complex are abstract simplicial complexes, that are defined on point cloud data in a metric space. And different $\varepsilon$ values will usually result in different $\check{C}$ech and Rips complexes. However, only $\check{C}ech$ complex preserves the homotopy information of the topological spaces formed by the $\varepsilon$-balls.

\subsubsection{Alpha Complex}

Alpha complex is used in computational geometry \citep{Edelsbrunner:1994,Alpha,Edelsbrunner:2010}. It is a set of subcomplexes from the Delaunay triangulation of a point set in $R^d$.

\begin{defn} Let X be a finite point set in $R^d$. Voronoi cell of a point x in X is the set of points $V_x \subseteq R^d$ for which x is the closest of the points in X, i.e. $V_x=\{u\in R^d: |u-x|\leq |u-x'|, \forall x' \in X\}$. 
\end{defn}

\begin{defn} The Delaunay complex of a finite set $X \in R^d$ is defined as the nerve of the Voronoi diagram.
\end{defn}

We define $R_{\varepsilon}(x)$ to be the intersection of the Voronoi cell $V_x$ with $B(X,\varepsilon)$, i.e.,
$R_{\varepsilon}(x)=V_x\bigcap B(X,\varepsilon)$ for evert $x \in X$.

\begin{defn}
The alpha complex $\mathcal{A}_{\varepsilon}(X)$ is definied as the nerve of the covers formed by $R_{\varepsilon}(x)$ for every $x \in X$, i.e. $\mathcal{A}_{\varepsilon}(X):=\left\{ \sigma \in X \mid \bigcap_{x \in \sigma}R_{\varepsilon}(x) \neq \emptyset \right\}$. 
\end{defn}
The geometric realization of alpha complex $\mathcal{A}_{\varepsilon}(X)$ is homotopy equivalent to the union of the closed balls $\bigcup_{x \in X} B(X,\varepsilon)$.

\subsubsection{Clique Complex}
Clique complex, also known as flag complex, is derived from graph and network models \citep{frohmader2008face,giusti2015clique,zomorodian2010tidy}. Its definition is very straightforward. For a
graph $G$, a complete graph is a graph, of which any two vertices are connected by an edge and a clique is a set of vertices that forms a complete subgraph in $G$. A $k$-clique complex is formed from a clique with $k+1$ vertices. Since subsets of cliques are also cliques, the clique complex is a simplicial complex.

\begin{figure}[!ht]
	\begin{center}
		\begin{tabular}{c}
			\includegraphics[width=0.5\textwidth]{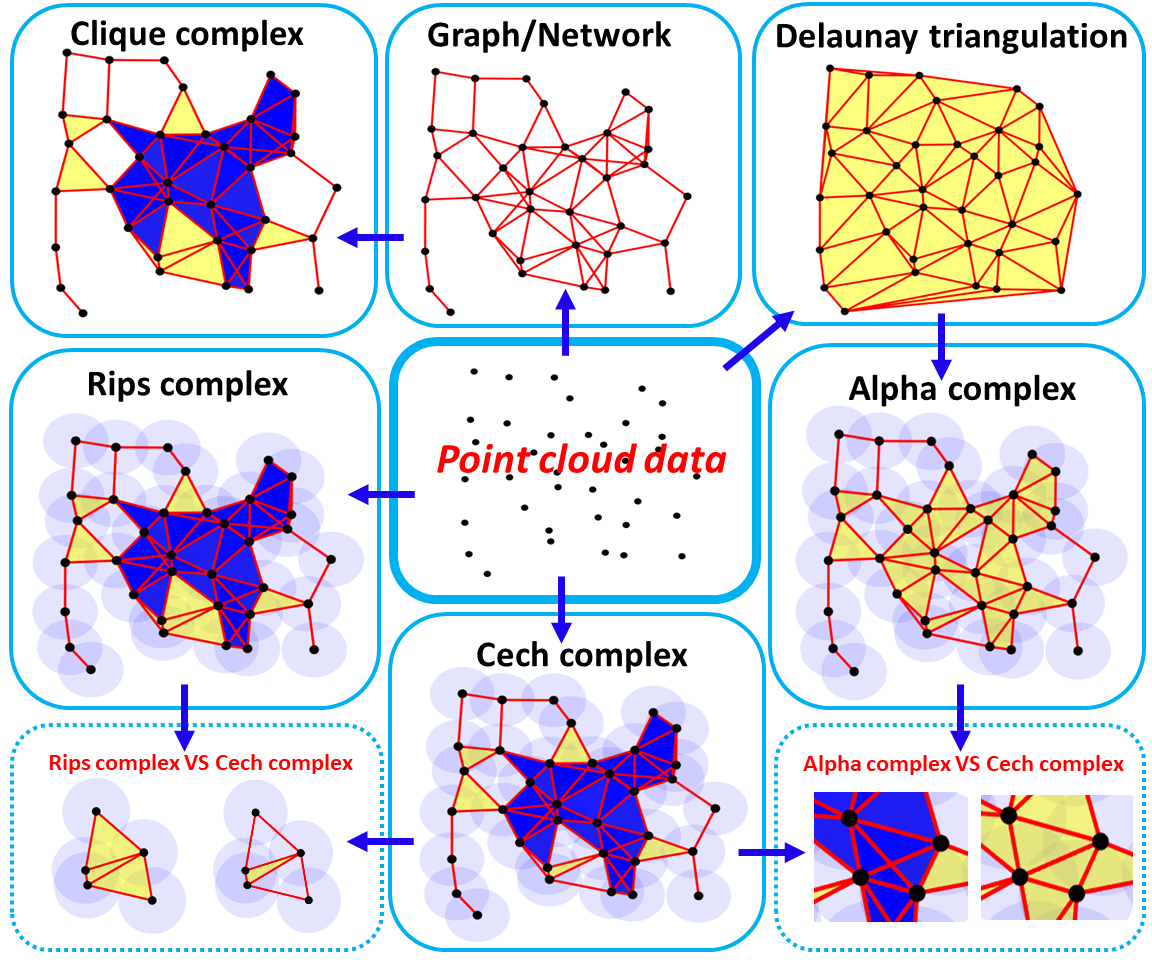}
		\end{tabular}
	\end{center}
	\caption{The comparison of $\check{C}$ech complex, Vietoris-Rips complex, Alpha complex and Clique complex. Essentially, Rips and $\check{C}$ech simplexes are different as Rips only requires the ``overlapping" of any two vertices to form a simplex. Clique complex is based on cliques that are derived from networks or graphs. Alpha complex is a subset of Delaunay triangulation and only suitable for subjects in Euclidian spaces. Interestingly, the same point cloud data, under different rules, can be used to generated different simplicial complexes, indicating that physical properties of data and its underly topological space is important for the construction of the simplicial complex.
	}
	\label{fig:simplicial_complex}
\end{figure}

Figure \ref{fig:simplicial_complex} illustrate the the comparison between $\check{C}$ech complex, Vietoris-Rips complex, Alpha complex and Clique complex. A point cloud data can be used to generate different simplicial complexes. The construction of the simplicial complex is highly dependent on the physical properties of data and the underlying topological space.

\subsubsection{Cubical Complex}
Image and volumetric data are arguably the important types of data in scientific computing. It is of key importance in computational topology to build up suitable simplicial complexes from these data. Since image and volumetric data is made from pixels and voxels, which are two and three dimensional cubes, it is natural to come up with the idea of cubical complex \citep{Kaczynski:2004}.

\begin{defn}
	An elementary interval is a closed interval $I \subset R$ of the form  $I=[l,l+1]$ or $I=[l,l]$, for some $l \in Z$. Here $[l,l]$ is an interval contains only one point. Elementary intervals that consist of a single point are degenerate.
\end{defn}

\begin{defn}
	An elementary cube $Q$ is a finite product of elementary intervals, that is $Q=I_1\times I_2\times \cdots \times I_d \subset R^d$, where each $I_i$ is an elementary interval. The set of all elementary cubes in $R^d$ is denoted by $\mathcal{K}^d$. The set of all elementary cubes is denoted by $\mathcal{K}$, $\mathcal{K}:=\bigcup_{d=1}\mathcal{K}^d$
\end{defn}

\begin{defn}
	A set $K \subset R^d$ is a cubical complex if $K$ can be written as a finite union of elementary cubes.
\end{defn}

\subsubsection{CW Complex}
CW complex plays an important role in topology and it has more freedom in the shape of the simplex. More specifically, we can define a closed unite ball in d-dimensional Euclidean space as $B^d=\{x \in R^d : |x| \leq 1\}$ and its boundary is the unit $(d-1)$-sphere $S^{d-1}=\{x \in R^d :|x|=1\}$. A $d$-cell is a topological space which is homemorphic to $B^d$, i.e., there exist a bijection mapping between $d$-cell and $B^d$. Roughly speaking, a cell can be viewed as a simplex with more complicated shape. The boundary of $d$-cell is a $(d-1)$-cell and a CW complex is formed when these cells satisfied the two basic conditions as in Definition \ref{def:two_conditions}.

\begin{defn}
	A finite CW complex is any topological space X such that there exists a finite nested sequence $\emptyset \subset X^0 \subseteq X^1 \subseteq \cdots \subseteq X^n=X$ such that for each $i=0, 1, 2,...,n$, $X_i$ is the results of attaching a cell to $X_{i-1}$. The sequence above is called the CW decomposition of X.
\end{defn}

\subsubsection{Morse-Smale Complex}

The Morse-Smale complex is a special kind of CW complex. Morse-Smale complex is defined usually on a manifold. It explores the manifold properties through a specially-designed function on it.

\begin{defn}
	Let f be a $C^2$ continuous function on a manifold $\mathbb{M}$. A point $p\in \mathbb{M} $ is a critical point if its derivatives with respect to a local coordinate system on $\mathbb{M}$ vanish.
\end{defn}

Morse theory studies the non-degenerate critical points, known as Morse points. Usually, the topological property of the critical points can be derived from their Hessian matrix,
\begin{eqnarray*}\label{eq:hessian}
H(f)=\left[ \begin{array}{cccc}
\frac{\partial f}{\partial x_1^2}   & \frac{\partial f}{\partial x_1 \partial x_2} &... &  \frac{\partial f}{\partial x_1 \partial x_n} \\
\frac{\partial f}{\partial x_2 \partial x_1}   & \frac{\partial f}{\partial x_2^2} &... &  \frac{\partial f}{\partial x_2 \partial x_n} \\	
\vdots  &\vdots  &\ddots &  \vdots \\
\frac{\partial f}{\partial x_n \partial x_1}   & \frac{\partial f}{\partial x_n \partial x_2} &... &  \frac{\partial f}{\partial x_n^2}
\end{array} \right].
\end{eqnarray*}
Let $\lambda_1 \leq \lambda_2 \leq ... \leq \lambda_n$ be the eigenvalues of $H(f)$. A critical point is degenerate when one of the eigenvalues equals to zero, otherwise it is non-degenerate. For non-degenerate critical point, its  number of negative eigenvalues is defined as the index of this critical point. This index is highly related to the topology of the manifold.

\begin{lem}
	(Morse Lemma) Let p be a non-degenerate critical point of $f$ with index $\lambda$ and let $c = f(p)$, there exists a local coordinate system
	$y = (y_1, y_2,..., y_n)$ in a neighborhood of $p$ with $p$ as its origin and $f(y)=c-y_1^2-y_2^2-...-y_{\lambda}^2+y_{\lambda+1}^2+y_{\lambda+2}^2+...+y_n^2$. 
\end{lem}

\begin{defn}
	Morse function is a smooth function on a manifold, such that all critical points are non-degenerate, and the critical points have distinct function values.
\end{defn}

Based on the Morse function, the gradient flow can be calculated and used to decompose a manifold. More specifically, at a point $p$ in $\mathbb{M}$, we can define its tangent plane as $T_p(\mathbb{M})$, $\gamma$ as any curve passing through it, and tangent vector as $v_p \in T_p(\mathbb{M})$ with $v$ the velocity of curve $\gamma$. We can explore how a Morse function changes in the direction specified by the tangent vector.

Mathematically, for a curve $\gamma(s)$ passing through p and tangent to $v_p \in T_p(\mathbb{M})$, the gradient $\nabla f$ of a Morse function $f$ along this curve is $\frac{d\gamma}{ds} \cdot \nabla f= \frac{d(f (\gamma))}{ds}$.

\begin{defn}
	An integral line $\gamma(s) : R \rightarrow \mathbb{M}$ is a maximal path whose tangent vectors agree with the gradient, i.e., $\frac{\partial}{\partial s} \gamma(s)=\nabla f(\gamma(s))$ for all $s \in R$. We call ${\rm org} \gamma = lim_{s \rightarrow -\infty} \gamma(s)$ the ${\rm origin}$ and ${\rm dest} \gamma = {\rm lim}_{s \rightarrow +\infty} \gamma(s)$ the ${\rm destination}$ of the path $\gamma$.
	Integral lines have the following properties. First, two integral lines are either disjoint or the same. Second, the integral lines cover all of $\mathbb{M}$. Third, the limits ${\rm org} p$ and ${\rm dest} p$ are critical points of $f$.
\end{defn}

\begin{defn}
	The stable manifold $S(p)$ and the unstable manifold $U(p)$ of a critical point p are defined as $S(p)=\{p\} \cup \{x \in \mathbb{M}~|~x \in {\rm im} \gamma,~{\rm dest} \gamma = p\}$ and $U(p)=\{p\} \cup \{x \in \mathbb{M}~|~x \in {\rm im} \gamma,~{\rm org} \gamma = p\}$, respectively, where $im \gamma$ is the image of the path $\gamma \in \mathbb{M}$.
\end{defn}
Both sets of manifolds decompose $\mathbb{M}$ into open cells.

\begin{defn}
	The stable manifold S(p) of a critical point p with index $i=i(p)$
	is an open cell of dimension dim $S(p) = i$.
\end{defn}

A Morse function is a Morse-Smale function if the stable and unstable manifolds intersect only transversally.

\begin{defn}
	(Morse-Smale complex) Connected components of sets $U(p_i)\bigcap S(p_j)$ for all critical points $p_i, p_j \in \mathbb{M}$ are Morse-Smale cells. We refer
	to the cells of dimension 0, 1, and 2 as vertices, arcs, and regions, respectively. The collection of Morse-Smale cells form a complex, the Morse-Smale
	complex.
\end{defn}

\subsection{Persistent Homology Analysis} \label{sec:PHA}

\subsubsection{Homology}

\begin{defn}
	An orientation of a $k$-simplex $\sigma^k=\{v_0, v_1, v_2, \cdots, v_{k+1} \}$ is an equivalence class of orderings of the vertices of $\sigma^k$, where $(v_0, v_1, v_2, \cdots, v_{k+1}) \sim (v_{\tau(0)}, v_{\tau(1)}, v_{\tau(2)}, \cdots, v_{\tau(k+1)})$	are equivalent orderings if the parity of the permutation $\tau$ is even. We denote an oriented simplex by $[\sigma^k]=[v_0, v_1, v_2, \cdots, v_{k+1} ]$.
\end{defn}

\begin{defn}
	A $k$-th chain group $C_k(K)$ of a simplicial complex $K$ is an Abelian group made from all $k$-chains from the simplicial complex $K$ together with addition operation.
\end{defn}

A linear combination of $k$-simplexes forms a $k$-chain $c$, i.e., $c=\sum_{i}\alpha_i\sigma^k_i$. In computational topology, we usually assume $\{ \alpha_i \in Z_2 \} $. All $k$-chains from the simplicial complex $K$ together with addition operation (modulo-2) will form an Abelian group $C_k(K, \mathbb{Z}_2)$. Homomorphism can be defined between these Abelian groups.

\begin{defn}
	The boundary operator $\partial_k: C_k \rightarrow C_{k-1}$ of an oriented $k$-simplex $[\sigma^k]=[v_0,v_1,v_2,\cdots,v_k]$ can be denoted as	$\partial_k [\sigma^k] = \sum^{k}_{i=0} [ v_0, v_1, v_2, \cdots, \hat{v_i}, \cdots, v_k ]$. Here $[v_0, v_1, v_2, \cdots ,\hat{v_i}, \cdots, v_k ]$ means a $(k-1)$ oriented simplex, which is generated by the original set of vertices $v_0,v_1,v_2,\cdots,v_{i-1},v_{i+1},\cdots, v_k$ except $v_i$.
\end{defn}
The boundary operator satisfies several properties, including $\partial_0= 0$ and $\partial_{k-1}\partial_k= 0$. With the boundary operation, we can define the $k$-th cycle group $Z_k$ and the $k$-th boundary group $B_k$ as
$$
Z_k={\rm Ker}~ \partial_k=\{c\in C_k \mid \partial_k c=0\}, \quad B_k={\rm Im} ~\partial_{k+1}= \{ c\in C_k \mid \exists d \in C_{k+1}: c=\partial_{k+1} d\}.
$$
Since $\partial_{k-1}\partial_k= 0$, cycle and boundary groups satisfy $B_k\subseteq Z_k$.

\begin{defn}
	The $k$-th homology group is the quotient group $H_k=Z_k/B_k$. The rank of $k$-th homology group $H_k$ is called $k$-th Betti number and satisfies $\beta_k = {\rm rank} ~H_k= {\rm rank }~ Z_k - {\rm rank}~ B_k$.
\end{defn}
Geometrically, we can regard $\beta_0$ as the number of isolated components, $\beta_1$ the number of one-dimensional loops, circles, or tunnels and $\beta_2$ the number of two-dimensional voids or holes.

\subsubsection{Persistent Homology}

\begin{figure}[!ht]
	\begin{center}
		\begin{tabular}{c}
			\includegraphics[width=0.9\textwidth]{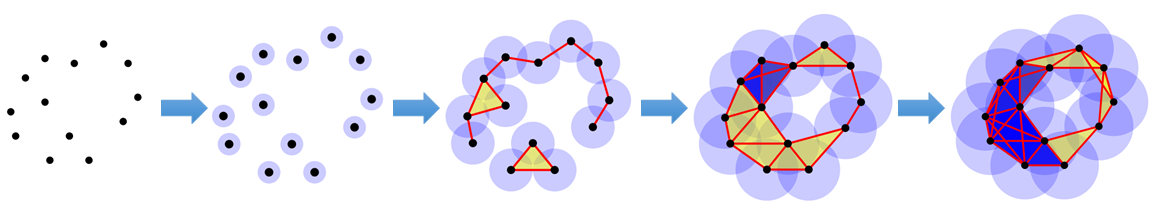}
		\end{tabular}
	\end{center}
	\caption{The illustration of a filtration process in PH. Each point is associated with a same-sized sphere. And the sphere radius is used as the filtration parameter. With the enlargement of the radius size, a series of nested simplicial complexes are generated.
	}
	\label{fig:filtration}
\end{figure}

The essence of PH is its filtration process, during which a series of topological spaces in different scales are generated (see Figure \ref{fig:filtration}).

\begin{defn}
	A filtration of a simplicial complex $K$ is a nested sequence of subcomplexes, $\emptyset \subset K^0 \subseteq K^1 \subseteq \cdots \subseteq K^n=K$. We call a simplicial  complex $K$ with a filtration a filtered complex.
\end{defn}

\begin{defn}
	Let $K^l$ be a filtration of a simplicial  complex $K$, and the $Z_k^l=Z_k(K^l)$ and $B_k^l=B_k(K^l)$ be the $k$-th cycle and boundary group of $K^l$, respectively. The $k$-th homology group of $K^l$ is $H^l_k=Z^l_k/B^l_k$. The $k$-th Betti number $\beta^l_k$ of $K^l$ is the rank of $H^l_k$.
\end{defn}

\begin{defn}
	The $p$-persistent $k$-th homology group at filtration time $l$ can be represented as
	$$
	H^{l,p}_k=Z^i_k/(B_k^{l+p}\bigcap Z^l_k).
	$$
	The $p$-persistent $k$-th Betti number $\beta_k^{l,p}$ of $K^l$ is the rank of $H^{l,p}_k$.
\end{defn}
Essentially, persistence gives a geometric measurement of the topological invariant.

\begin{figure}[!ht]
	\begin{center}
		\begin{tabular}{c}
			\includegraphics[width=0.9\textwidth]{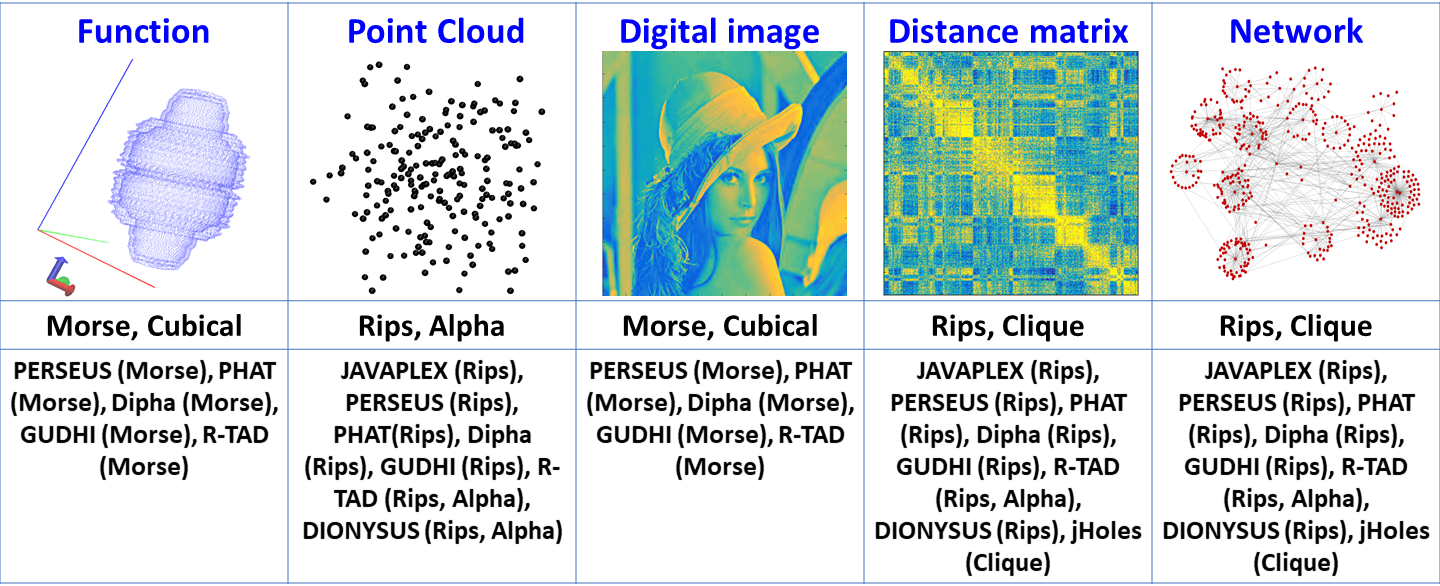}
		\end{tabular}
	\end{center}
	\caption{The illustration of the simplicial complexes and PH softwares for different types of data. Depending on the data types, different simplicial complexes should be used to capture the intrinsic topology of the data. For instance, Rips complex can be used for point cloud data. With a suitable filtration parameter, PH can be further employed to extract the multiscale topological information. Various softwares are available for PH analysis.
	}
	\label{fig:data_SimCom}
\end{figure}

\paragraph{Softwares for PH}\label{Overview of available software packages for PH}
Various softwares, including JavaPlex \citep{javaPlex}, Perseus  \citep{Perseus}, Dipha \citep{Dipha}, Dionysus \citep{Dionysus}, jHoles \citep{Binchi:2014jholes}, GUDHI \citep{gudhi:FilteredComplexes}, Ripser \citep{bauer2017ripser}, PHAT \citep{bauer2014phat}, DIPHA \citep{bauer2014distributed}, R-TDA package \citep{fasy:2014introduction}, etc, have been developed. These softwares are based on different types of simplicial complexes and can be used to process different types of data (see Figure \ref{fig:data_SimCom}).

\begin{figure}[!ht]
	\begin{center}
		\begin{tabular}{c}
			\includegraphics[width=0.9\textwidth]{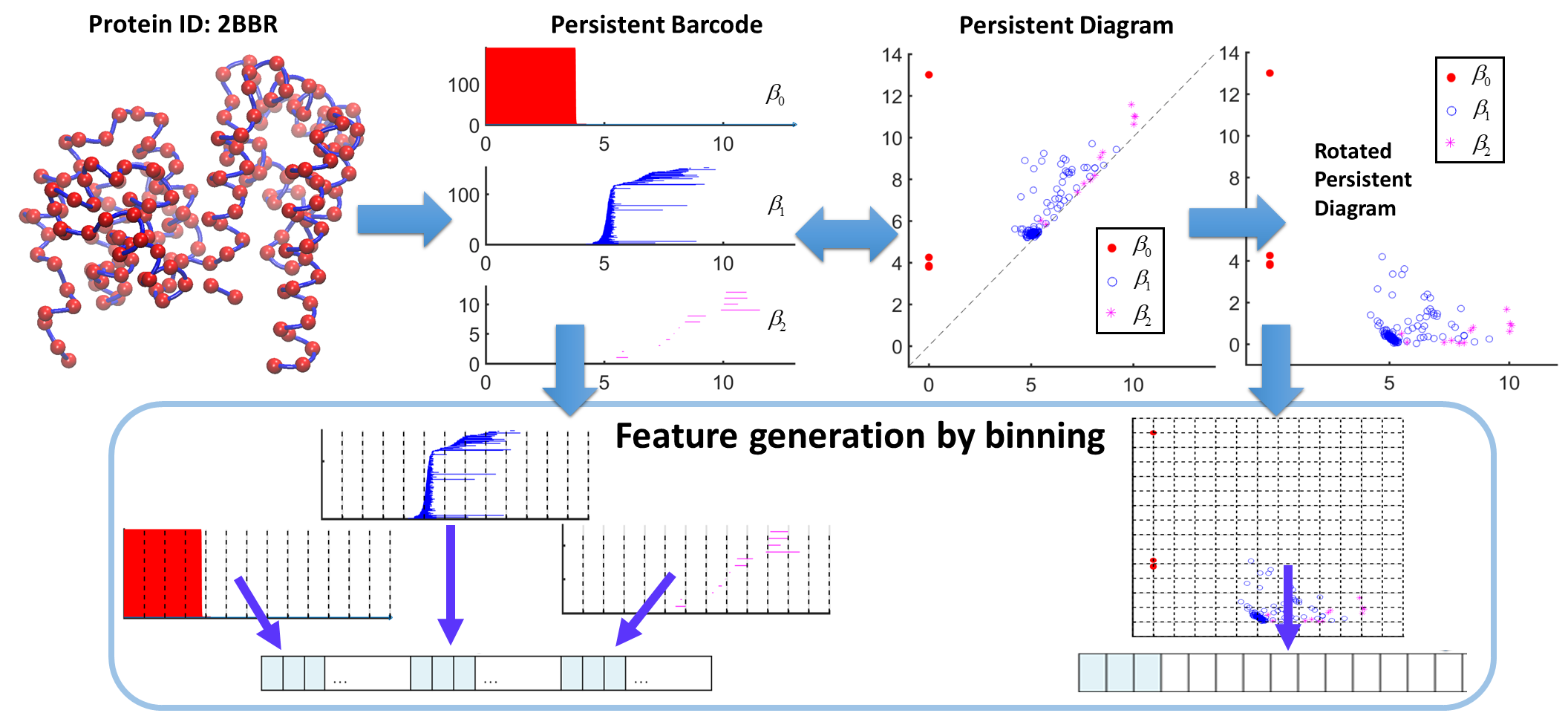}
		\end{tabular}
	\end{center}
	\caption{The illustration of the topological feature generation by binning approach. Essentially, binning is to discretize the PBs/PDs into equally-sized elements, count the corresponding Betti numbers (as in PB) or PD points, and concatenate them into a feature vector. The feature vector can then be used in machine learning models.
	}
	\label{fig:barcodes_binning}
\end{figure}

\paragraph{Persistent barcodes and persistent diagram}\label{Barcode representation}
The results from PH can be represented as pairs of birth times (BTs) and death time (DTs). More specifically, for each generator (topological invariant), BT and DT represent the filtration value at which the generator is born and vanished, respectively. They always come in pairs and can be denoted as
$$
l^k_j=\{a^k_j, b^k_j\} ~{\rm with }~  k \in \mathbb{N} ~{ \rm and}~ j \in \{1, 2, ..., N_k\}.
$$
We use $a^k_j$, $b^k_j$, $l^k_j$ to represent BT, DT, persistence for the $j$-th topological generator of the $k$-th dimensional Betti number, respectively. Parameter $N_k$ is the total number of $k$-th dimensional topological generators. Moreover, dimension number $k$ is always a nonnegative integer (i.e., $k \in \mathbb{N}$). These notations will be consistently used in the following sections unless stated otherwise. To avoid repetition, we will omit the definition of domain for the dimension parameter $k$ and always assume it to be $k \in \mathbb{N}$. For simplification, we also define the set of $k$-th dimensional barcodes as,
$$
L_{k}= \{ l^k_j=\{a^k_j, b^k_j\}_{j \in \{1, 2, ..., N_k\}}\}.
$$
The topological persistence is usually represented in two ways, i.e., persistent barcode (PB) \citep{Ghrist:2008} and persistent diagram (PD) \citep{Mischaikow:2013}. In PB, each $l^k_j$ is treated as a bar, i.e., the line $[l^k_j]=[a^k_j, b^k_j]$. In PD, each $l^k_j$ is considered as a 2D point with coordinate $ \vec{l^k_j}=(a^k_j, b^k_j)$. The upper panel of Figure \ref{fig:barcodes_binning} illustrates the PB and PD for a protein structure with PDB ID: 2BBR.  The rotated PD is also widely used. The essential idea is to employ a rotational function  $T(x,y)=(x, y-x)$ on the data points in PD. In this way, a point $(a^k_j, b^k_j)$ in PD $(a^k_j, b^k_j)$, becomes $T(\vec{l^k_j})=(a^k_j, b^k_j-a^k_j)$ after rotation. The rotated PD is denoted as
$$
T(L_k)=\{ T(\vec{l^k_j})= (a^k_j, b^k_j-a^k_j)_{ j \in \{1, 2, ..., N_k\}}\}.
$$

\paragraph{Persistent local homology}
A closely related model, the persistent local homology \citep{bendich2007inferring,bendich2012local,ahmed2014local,bendich2015multi,fasy2016exploring}, is proposed based on the algebraic topological concept called local homology groups \citep{munkres2018elements}. Generally speaking, local homology is used to evaluate the local structure of a topological space. If the topological space is topological manifold, the local homology groups for each point are the same. More specifically, we can let $X$ be a topological space embedded in $R^d$. For a fixed point $p \in R^d$, we can define a ball center at $p$ with radius $\epsilon$ $B(p, \epsilon)=\{x \in R^d : |x| \leq 1\}$. The $k$-th local homology group of $X$ with center point $p$ and radius $\epsilon$ is denoted as $LH_k(X, p, \epsilon)=H_k(X\bigcap B(p, \epsilon))$. Persistent local homology can be used in dimension reduction and manifold dimension detection.

\subsubsection{Persistent Functions}

Based on PB/PD results, different functions are proposed to represent or analyze the topological information \citep{Carlsson:2009,bubenik:2015,Chintakunta:2015}. They can be treated as functional data for subsequent analysis \citep{Ramsay1997,Horvath2012} or used to deduce the structured features as shown in Section \ref{sec:bin}-\ref{sec:2d3d} below.
\begin{defn}
	The persistent Betti number (PBN) or Betti curve is defined as the summation of all the $k$-th dimensional barcodes,
	$$
	f(x;L_{k})= \sum_{j} \chi_{[a^k_j,b^k_j]}(x)
	$$
	Function $\chi_{[a^k_j, b^k_j]}(x)$ is a step function, which equals to one in the region $[a^k_j, b^k_j]$ and zero otherwise.
\end{defn}
This equation transforms all the $k$-dimensional barcodes into a one-dimensional function. However, the persistent Betti number is not continuous.

\begin{defn}
	A continuous persistent Betti function is proposed as follows,
	$$\label{eq:PBF}
	f(x;L_{k})= \sum_{j} e^{-\left(\frac{x- \frac{a^k_j+b^k_j}{2}}{w_j (b^k_j-a^k_j)}\right)^2},
	$$
	where $w_j$ are weight values, which can be adjusted based on different purposes \citep{Xia:2017similarity}.
\end{defn}

Another way to represent the barcode result is persistence landscapes \citep{bubenik:2015}.
\begin{defn}
	A piecewise linear function can be defined on each individual barcode $\{ l^{k}_{j} \}$,
	\begin{eqnarray*}\label{eq:couple_matrix25}
	f(x,l^{k}_{j})=\begin{cases} \begin{array}{ll}
	0 &  {\rm if} \quad x \not\in (a^k_j, b^k_j);\\
	x-a^k_j  & {\rm if} \quad x \in (a^k_j, \frac{a^k_j+b^k_j}{2} ];\\		
	-x+b^k_j & {\rm if} \quad x \in [\frac{a^k_j+b^k_j}{2}, b^k_j).
	\end{array}
	\end{cases}
	\end{eqnarray*}
	The persistence landscape of $k$-dimensional barcode $L_k$ is the sequence functions $\lambda_m: R \rightarrow [0, \infty], m=1, 2, 3,...$, where $\lambda_m(x)$ is the $m$-th largest value of $\{f(x,l^{k}_{j})\}_{j=1}^{N_k}$. We set $\lambda_m(x)=0$ if the $m$-th largest value does not exist.
\end{defn}

Recently, persistent image has been proposed to represent PDs \citep{adams2017persistence}. The persistent image is derived from persistent surface function.
\begin{defn}
	The persistent surface for $k$-th dimensional barcodes $L_{k}$ is defined as
	$$
	\rho (x,y,L_{k})=\sum_j w(x,y,t_1,t_2)\phi(x,y,l^k_j),
	$$
	where the differentiable probability distribution function $\phi(x,y,l^k_j)=\frac{1}{2\pi \sigma^2}e^{-((x-a^k_j)^2+(y-(b^k_j-a^k_j))^2)/2\sigma^2}$ and the weight function with two constants $t_1$ and $t_2$ is
	\begin{eqnarray*}\label{eq:weight_function}
	w(x,y,t_1,t_2)=\begin{cases} \begin{array}{ll}
	0 &  {\rm if} \quad y \leq t_1;\\
	\frac{y-t_1}{t_2-t_1}  & {\rm if} \quad t_1<y<t_2;\\		
	1 & {\rm if} \quad y\geq t_2.
	\end{array}
	\end{cases}
	\end{eqnarray*}
	The persistent surface can be discretized by a Cartesian grid into image data. By the integration of persistent surface function over each grid (or pixel), a persistent image function is obtained.
\end{defn}

To measure the disorder of a system, persistent entropy has been proposed \citep{Merelli:2015topological,Chintakunta:2015,Rucco:2016,Xia:2018multiscale}.
\begin{defn}
	Persistent entropy is defined as
	$$S_k=\sum_j - p^k_j ln(p^k_j)$$
	with the probability function $p^k_j=\frac{b^k_j-a^k_j}{\sum_j (b^k_j-a^k_j)}$.
\end{defn}

\section{Persistent-Homology-based Machine Learning Methods}\label{Sec:PH_SL}

\subsection{Feature Extraction from PBs/PDs} \label{sec:feature}

The topological information within PBs/PDs needs to be converted to (structured) features, so that they can be inputted in machine learning algorithms. These features are usually represented as a quantitative vector. In this section, we summarize the common approaches of feature extraction from PBs/PDs.

\subsubsection{Barcode statistics}
The simplest way of feature generation is to collect the statistic measurement of the barcode properties, such as birth times (BTs), death times (DTs) and persistent lengths (PLs). These statistic measurements can be maximal, minimal, mean, variance, summation, among others. In the subsequent experiments, we refer to the 13 barcode statistics in  \citep{ZXCang:2015}, which are used to characterize topological structure information. For instance, first/second/third longest Betti-0/Betti-1/Betti-2 PLs, BT/DT for the longest Betti-1/Betti-2; number of Betti-1 bars located at filtration distance [4.5, 5.5] \AA, etc.

\subsubsection{Algebraic functions and tropical functions}
Algebraic combinations were also considered \citep{adcock2016ring}. For instance,
$\sum_j a_j^k(b_j^k-a_j^k)/N_k$, $(\max_{j}\{b_j^k\}-b_j^k)(b_j^k-a_j^k)/N_k$, $(a_j^k)^2(b_j^k-a_j^k)^4$, $( \max_{j}\{b_j^k\}-b_j^k)^2(b_j^k-a_j^k)^4/N_k$.
The average over each bar to eliminate the effects of large variations.

An alternative is to consider the min-plus and max-plus type coordinates \citep{kalivsnik2018tropical}. For instance,
$\max_j \{(b_j^k-a_j^k)\}$, $\max_{i<j}\{(b_i^k-a_i^k)+(b_j^k-a_j^k)\}$,  $\max_{i<j<m} \{(b_i^k-a_i^k)+(b_j^k-a_j^k)+(b_m^k-a_m^k)\}$, $\max_{i<j<m<n} \{(b_i^k-a_i^k)+(b_j^k-a_j^k)+(b_m^k-a_m^k)+(b_n^k-a_n^k)\}$, $\sum_j\{ (b_j^k-a_j^k)\}$.

\subsubsection{Binning approach} \label{sec:bin}
Structured feature vectors can be constructed by simple yet effective approach of binning. The idea is to discretize the filtration domain into various same size bins (see Figure \ref{fig:barcodes_binning}). For PBs, one dimensional bins are used, i.e., $[x_i, x_{i+1}], i=0, 1, ..., n$ with $x_0=0$ and $x_n=r_f$ (with $r_f$ the filtration ending value).

\begin{enumerate}\small
	\item For persistent Betti number and persistent Betti function, the binning process can be done easily by collecting their values at the grid points \citep{cang:2017topologynet,cang:2018representability}, that is, the set $\{f(x_i,L_k) |i=0, 1, ... , n; k=0,1,2\}$. In this way, we can systematically obtain $n$ features for each Betti number.
	
	\item For persistent landscapes $\lambda_m(x)$, the same binning process as stated above can be done repeatedly for $m$ times. The $m$ set of $\{\lambda_m(x_i) |i=0, 1, ... ,n \}$ values are used as features \citep{bubenik2017persistence}.
	
	\item For persistent image \citep{adams2017persistence}, the binning process is based on the rotated PD. Essentially, the computational domain of the rotated PD can be discretized in $n\cdot n$ grids (or pixels) and persistent image values can be evaluated on each grid. Therefore, there are totally $n \cdot n$ features in persistent image vector.
	
	\item The distribution functions of BTs, DTs and PLs can be discretized and used as feature vectors \citep{cang:2017topologynet,cang:2018representability}. The idea is very straightforward. For each interval  $[x_i, x_{i+1}]$, we can count the numbers of the $k$-th dimensional BTs, DTs, PLs located in this range, and denoted them as $N_{BT}^{k,i}$,$N_{DT}^{k,i}$,$N_{PL}^{k,i}$, respectively. The sets of count numbers $\{N_{BT}^{k,i}, N_{DT}^{k,i}, N_{PL}^{k,i} | i=0, 1, ..., n; k=0, 1, 2 \}$ can be assembled as a feature vector. It should be noticed that normally for Betti-0, the BTs are 0, thus DTs are usually equal to PLs. So only the set of $\{N_{PL}^{0,i} | i=0, 1, ..., n\}$ is used.
	
	\item Similarly, if we consider the PD representation \citep{cang:2017topologynet,cang:2018representability}, we can discretize the diagram domain into a mesh and count the number of points in each grid. We can also consider the rotated PD and count the point numbers in each grid. These numbers form a same-size vector and can be used as the feature vector.
\end{enumerate}

\subsubsection{Persistent codebook}
Instead of using the regular mesh, persistent codebooks are proposed by using different clustering or bagging of the PD points \citep{zielinski2018persistence,bonis2016persistence}.
\begin{enumerate}\small
	\item  The idea of persistent bag of words (P-BoW) is to assign the points in the rotated PD $T(L_k)=\{ T(\vec{l^k_j})= (a^k_j, b^k_j-a^k_j) \}$ to a precomputed codebook. In P-BoW, the K-means clustering algorithm is employed to cluster the points $T(L_k)$ into $n$ clusters, i.e., $NN(T(\vec{l^k_{i,m}}))=i (i=1,2,..,n; m=1,2,....N_i)$. Here $NN(x,y)=i$ means point $(x,y)$ belongs to cluster $i$ and $N_i$ is the total number of points in $i$-th cluster. Further we denote $\vec{x_i}=(x_i,y_i)$ as the cluster center and the P-BoW representation is $\vec{v}=(v_i=|\vec{x_i}|)_{i=1, 2, ..., n }$. To further consider the persistence information, the weight function $w(x,y,t_1,t_2)$ as in persistent surface is considered.
Usually, the parameters $t_1$ and $t_2$ are set to the persistence values corresponding to 0.05 quantile of birth time and 0.95 quantile of death time, respectively. The weighted K-means clustering provides a more adaptive codebook.
	
	\item  Persistent vector of locally aggregated descriptors (P-VLAD) uses the K-means clustering as in P-BoW and then the accumulated distance between a point $T(\vec{l^k_{i,j}})$ to the nearest codeword $\vec{x_i}$ is computed as following $\vec{v}_i=\sum_{ m=1,2,...,N_i} (T(\vec{l^k_{i,m}})-\vec{x_i})$ with $N_i$ the total number of points in $i$-th cluster. These $n$ vectors can be concatenated into a $2n$ dimensional vector. It can capture more information than the P-BoW.
	
	\item Persistent Fisher vector (PFV) is to characterize the rotated PD with a gradient vector from a probability model. We let $\lambda_{\rm GMM}=\{ w_i, \vec{\mu_i}, \sum_i; i=1,2,...,m \}$ be the set of parameters of Gaussian mixture models (GMMs), where $w_i$, $\mu_i$ and $\sum_i$ denote the weight, Guassian center and covariance matrix of $i$-th Gaussian, respectively. The likelihood that point $T(\vec{l^k_j})$ is generated by $i$-th Gaussian is denoted as $p_i(T(\vec{l^k_j})|\lambda_{\rm GMM})$,
	$$
	p_i(T(\vec{l^k_j})|\lambda_{\rm GMM})=\frac{exp\{-\frac{1}{2}(T(\vec{l^k_j})-\vec{\mu_i})'\sum_i^{-1}(T(\vec{l^k_j})-\vec{\mu_i})\}}{2\pi|\sum_i|^{\frac{1}{2}}}
	$$
and the function is
	$$
	\mathcal{F}(T(L_k)|\lambda_{\rm GMM})= \sum_{j}log\Big( \sum_i^N w_i p_i(T(\vec{l^k_j})|\lambda_{\rm GMM})\Big).
	$$
	When the covariance matrices are diagonalized, the derivations $\frac{\partial \mathcal{F}(T(L_k)|\lambda_{\rm GMM})}{\partial \vec{\mu_i}}$ and $\frac{\partial \mathcal{F}(T(L_k)|\lambda_{\rm GMM})}{\partial \sigma_i}$ can be evaluated. Here the $\sigma_i$ is the diagonal element of $\sum_i$. The PFV is the concatenation of these partial derivatives.
\end{enumerate}

\subsubsection{Persistent paths and signature features}
Recently, a new feature map for barcodes is proposed \citep{chevyrev2018persistence}. This model includes two major steps. First, to embed the barcode information into persistent path. Several embeddings are proposed, including landscapes embedding, envelope embedding, Betti embedding and Euler embedding. The persistent landscapes $\lambda_m(x)$ is a typical example of embbing barcodes into path information.  Second, from persistent path to tensor series. And tensor series can be used as feature vectors.

\subsubsection{2D/3D representation} \label{sec:2d3d}
Two or three dimensional matrixes can be generated from PD, multi-dimensional PH and element-specific PH \citep{cang:2017topologynet,cang:2018representability}. The persistent surface is an example of 2D representation. In multi-dimensional PH, we stack each persistent Betti number along the other filtration parameter, this results in a 2D image. If three filtration parameters are considered, a 3D matrix will be generated. Similarly, in element-specific PH, each type of atom combination will contribute seven distributions for BTs, DTs and PLs, all types atoms together form 7 matrixes. Unlike the feature vector representations, these matrix representation can be treated as images and combined directly with convolution neural network (CNN).

\subsubsection{Network layer for topological signature}
A parameterized network layer for topological signatures is proposed in \citep{hofer2017deep}. We denote $R^+=(0, \infty)$ and $R^+_0=[0, \infty)$. We let $\vec{\mu}=(\mu_0, \mu_1)^T \in R \times R_0^+$ and $\vec{\sigma}=(\sigma_0, \sigma_1) \in R^+ \times R^+$ and $\upsilon \in R^+$. We can define,
\begin{eqnarray*}
s(\vec{l^{k}_{j}},\vec{\mu}, \vec{\sigma}, \upsilon)=\begin{cases} \begin{array}{ll}
e^{-\sigma_0^2(a_j^k-\mu_0)^2-\sigma_1^2(b_j^k-\mu_1)^2},  & b_j^k \in[\upsilon, \infty);\\
e^{-\sigma_0^2(a_j^k-\mu_0)^2-\sigma_1^2( \ln(\frac{b_j^k}{\upsilon})+\upsilon-\mu_1)^2}, & b_j^k \in(0, \upsilon);\\		
0, & b_j^k=0.
\end{array}
\end{cases}
\end{eqnarray*}
and  $S(\vec{\mu}, \vec{\sigma}, \upsilon)=\sum_{j} s(\vec{l^{k}_{j}},\vec{\mu}, \vec{\sigma}, \upsilon)$. For a set of $N$ different parameter pairs $(\vec{\mu_i}, \vec{\sigma_i}), i=1,2,...,N$, we can concatenate all the corresponding function values together as a vector $ (S(\vec{\mu_i}, \vec{\sigma_i}, \upsilon))_{i=1,2,...,N} $. This vector is then the output of the network layer.

\subsection{Persistent-Homology-based Kernel Methods} \label{sec:kernels}
With the significant role of kernels in machine learning models, one can propose suitable kernel functions for the topological data to replace the role of features. Since similarity measures can be modified into kernels \citep{Chen2009}, we are keen to develop similarity measurement for PH. Moreover, some kernels are separately proposed from the topological perspectives.

\subsubsection{Similarity measurement}
The unique format of PH outcomes poses a great challenge for a meaningful metrics. To solve this problem, distance measurements or metrics \citep{mileyko2011probability,turner2014frechet,chazal2017robust,carriere2018metric,anirudh2016riemannian} have been considered, including Gromov-Hausdorff distance, Wasserstein distance \citep{mileyko2011probability}, bottleneck distance \citep{mileyko2011probability}, probability-measure-based distance \citep{chazal2011geometric,marchese2017signal,chazal2017robust}, Fisher information metric \citep{anirudh2016riemannian}.

\begin{defn}
	The bottleneck distance between two sets of barcodes $L_{k}$ and $L_{k}'$ is defined as
	$$
	d_B(L_{k},L_{k}')= \underset{\gamma}{\inf} \, {\underset{j}{\sup} \| l^{k}_{j}-\gamma(l^{k}_{j}) \|_{\infty}},
	$$
	here $\gamma$ ranges over all bijections from barcodes $L_{k}$ to barcodes $L_{k}'$. The distance between two barcodes $L_{k,j}$ and $L_{k,j}'$ is defined as $\| l^{k}_{j}-l^{k}_{j'} \|_{\infty}=\max \left\{ |a_{k,j}-a_{k,j'}|, |b_{k,j}-b_{k,j'}| \right\}$.
\end{defn}

\begin{defn}
	The p-Wasserstein distance between two barcodes is defined as follows,
	$$
	d_{W,p}(L_{k},L_{k}')=\underset{\gamma}{\inf} \left[ \sum_{j} \|l^{k}_{j}-\gamma(l^{k}_{j}) \|_{\infty}^p \right]^{\frac{1}{p}}.
	$$
	Again $\gamma$ is bijection from barcodes $L_{k}$ to barcodes $L_{k}'$. Parameter $p$ is a positive integer.
\end{defn}

\begin{defn}
	The Sliced Wasserstein (SW) distance between two barcodes is defined as follows. For a $\theta \in [-\pi/2, \pi/2]\}$, we can define a line $f(\theta)= \{\lambda (\cos \theta,\sin \theta) | \lambda \in R\}$. Further $\pi_\theta : R^2 \rightarrow f(\theta)$ is defined as the orthogonal projection of a point onto this line and the $\pi_\Delta$ is the orthogonal projection onto the diagonal line (i.e., $\theta=\pi/4$). If we denote $\mu(\theta,L_{k}) =\sum_{j} \delta_{\pi_\theta(\vec{l^k_j})} $ and $\mu_{\Delta}(\theta, L_{k})=\sum_{j} \delta_{\pi_\theta \circ \pi_\Delta(\vec{l^k_j})}$, the SW distance is defined as,
	$$
	SW(L_{k},L_{k}')= \frac{1}{2\pi}\int \mathcal{W}(\mu(\theta,L_{k})+\mu_{\Delta}(\theta,L_{k}'),\mu(\theta,L_{k}')+\mu_{\Delta}(\theta,L_{k}) )d\theta,
	$$
	where the function $\mathcal{W}(x,y)$ is the generic Kantorovich formulation of optimal transport.
\end{defn}

\subsubsection{Topological kernel}
From the topological perspectives, various PH-based kernels have also been proposed, which are listed as follows.


\begin{defn}
	Persistent scale space kernel (PSSK) is proposed by \citep{reininghaus2015stable} as follows.
	$$
	\kappa_{\rm PSSK}(L_{k},L_{k}',\sigma)=\frac{1}{8\pi \sigma}\sum_{\vec{l^{k}_{j}}\in L_{k},~\vec{l^{k}_{j'}}\in L_{k}'}  e^{-\frac{\| \vec{l^{k}_{j}}-\vec{l^{k}_{j'}}\|^2}{8\sigma}}-e^{-\frac{\| \vec{l^{k}_{j}}-\overline{\vec{l^{k}_{j'}}}\|^2}{8\sigma}},
	$$
	where $\overline{\vec{l^{k}_{j'}}}$ is $\vec{l^{k}_{j'}}$ mirrored at the diagonal.
\end{defn}

\begin{defn}
	Universal persistence scale-space kernel (u-PSSK) is proposed by \citep{kwitt2015statistical} as follows.
	$$
	\kappa_{\rm uPSSK}(L_{k},L_{k}',\sigma)=e^{\kappa_{\rm PSSK}(L_{k},L_{k}',\sigma)}
	$$
\end{defn}

\begin{defn}
	Persistent weighted Gaussian kernel (PWGK) is proposed by \citep{kusano2016persistence} as follows.
	$$
	\kappa_{\rm PWGK}(L_{k},L_{k}',\sigma)=\sum_{\vec{l^{k}_{j}}\in L_{k},~\vec{l^{k}_{j'}}\in L_{k}'} w_{arc}(\vec{l^{k}_{j}})w_{arc}(\vec{l^{k}_{j'}})e^{-\frac{\|\vec{l^{k}_{j}}-\vec{l^{k}_{j'}}\|}{2\sigma^2}}
	$$
	here $w_{arc}(\vec{l^{k}_{j}})=\arctan(C (b^{k}_{j}-a^{k}_{j})^p)$ with parameter $C$ and $p$ all positive value.
\end{defn}

\begin{defn}
	Geodesic topological kernel (GTK) is proposed by \citep{padellini2017supervised} as follows.
	$$
	\kappa_{\rm GTK}(L_{k},L_{k}',\sigma)=e^{\left\{\frac{1}{h} d_{W,2}(L_{k},L_{k}')^2 \right\}}
	$$
	with $h>0$ and $d_{W,2}(L_{k},L_{k}')$ the 2-Wasserstein distance. Similarly, the geodesic Laplacian kernel (GLK) is defined as
	$$
	\kappa_{\rm GLK}(L_{k},L_{k}',\sigma)=e^{\left\{\frac{1}{h} d_{W,2}(L_{k},L_{k}') \right\}}
	$$
\end{defn}

\begin{defn}
	Persistence Fisher Kernel is proposed by \citep{le2018riemannian} as follows. We let
	$$
	\rho (x,y,L_{k})=\frac{1}{Z}\sum_j^{N_k}N(x,y|\vec{l^{k}_{j}},\sigma).
	$$
	with $Z= \int \sum_j^{N_k}N(x,y|\vec{l^{k}_{j}},\sigma) dxdy$ and $N(x,y|\vec{l^{k}_{j}},\sigma)=\frac{1}{\sqrt{2\pi \sigma^2}}e^{-\frac{(x-a^k_j)^2+(y-b^k_j)^2}{2\sigma^2}}$ is the normal distribution. Further we define the Fisher information metric (FIM) between $L_{k}$ and $L_{k}'$  as
	$$
	d_{\rm FIM}(\rho (x,y,L_{k}),\rho (x,y,L_{k}'))=\arccos \Big( \int \sqrt{\rho(x,y,L_{k} \cup \pi_\Delta (L_{k}')) \rho(x,y,L_{k}' \cup \pi_\Delta (L_{k}))} dxdy \Big)
	$$
	The persistent Fisher kernel (PFK) is defined as,
	$$
	\kappa_{\rm PFK}(L_{k},L_{k}')= e^{-t_0 d_{\rm FIM}(\rho (x,y,L_{k}),\rho (x,y,L_{k}'))}.
	$$
	Parameter $t_0$ is a positive scale value. The PFK is positive definite.
\end{defn}

\begin{defn}
	Three persistent-landscape-based kernels are defined as follows.
	\begin{enumerate}
		\item Global PH kernel (GPHK): For two persistent landscape functions $\lambda_m(x)$ and $\lambda'_m(x)$, the GPHK is defined as
		$$
		\kappa_{\rm GPHK}(\lambda_m,\lambda'_m)=\langle\lambda_m, \lambda'_m\rangle=\int \lambda_m(x)\lambda'_m(x)dx.
		$$
		\item Multiresolution PH kernel (MPHK): Let $B(c,r)$ be the ball center at point $c$ with radius $r$. We define a local homology by only considering the points within the ball $B(c,r)$ and filtration process is only performed for these special points. Local persistent landscape functions $\lambda_B(c,r)$ can be derived. If we choose the ball centers from a uniform distribution ($\mathcal{P}_C$) over the whole point cloud data $X$, the expected persistent landscape is obtained by $\pi_r=\mathbb{E}_{c\in \mathcal{P}_C} \lambda_{B(c,r)}$. To obtain topological information of $X$ at finer resolutions, we consider $M$ decreasing radius $r_1 >r_2 > \cdots > r_M $, and compute the corresponding expected persistent functions $\pi_{r_1}, \pi_{r_2},..., \pi_{r_M}$. We can define the MPHK as
		$$
		\kappa_{\rm MPHK}= \sum_i^M w_i \langle \pi_{r_i},\pi'_{r_i} \rangle,
		$$
		where the parametes $w_i$ are nonnegative weights to combine different resolutions. For instance, they can be chosen as $w_i=(r_1/r_i)^3$.
		\item Stochastic multiresolution PH kernel (SMURPHK): It aims to reduce the computational cost of the MPHK above. Essentially, $\pi_r=\mathbb{E}_{c\in \mathcal{P}_C} \lambda_{B(c,r)}$ can be approximated by considering $n$ centers in the probability function $\mathcal{P}_C$ and $s$ bootstrap samples in $\lambda_{B(c,r)}$. In this way, an average persistent landscape function can be defined as $\bar{\lambda}_r=\frac{1}{ns}\sum^n_{i=1} \sum^s_{j=1} \lambda_{b_j}(c,r)$. Then the SMURPHK is defined as
		$$
		\kappa_{\rm SMURPHK}=\sum_i^M w_i \langle \bar{\lambda}_{r_i},\bar{\lambda'}_{r_i} \rangle.
		$$
	\end{enumerate}
\end{defn}

\subsection{Machine Learning Algorithms} \label{sec:MLalgo}

\subsubsection{Support Vector Machine (SVM)}\label{svm}
Support vector machines (SVM), proposed in  \citep{Cortes:1995}, is a supervised learning algorithm that finds the optimal hyperplane/decision boundary of the form $\mathbf{w^Tx}+b=0$ in the context of 2-class classification problems. This creates a quadratic optimisation problem to solve for coefficients $\mathbf{w}$ and $b$. The primal problem is formulated as:
$$
\underset{\mathbf{w},b}\min \quad \frac{1}{2}{\mathbf{\left |w \right |}}^2 +C\sum\xi(\mathbf{w};\mathbf{x_i},y_i),
$$
where $\xi(\mathbf{w};\mathbf{x_i},y_i)$ is a convex non-negative loss function. The above L2-regularised classifier can also be modified into an L1-regularised classifier together with different loss function used.
Rewriting the problem to involve Lagrange multipliers converts the primal problem that optimises (minimises) with respect to coefficients $\mathbf{w}$ into the dual problem that optimises (maximises) with respect to the multipliers.
Detailed derivation of the optimal hyperplane algorithms and soft margin classifiers can be found in \citep{Cortes:1995}. The method is popular for its robustness to individual observations, compared to discriminant analysis (with the use of support vectors, namely only the training observations that lie within the margins are used).

In some cases, however, if the data is linearly inseparable, it requires the use of feature expansion, i.e. to map the data into a higher dimensional space $\Psi:\mathbb{R^\mathcal{N}}\to{\mathbb{R^\mathcal{M}}}, \mathcal{N}<\mathcal{M}$. Such a transformation involves the use of kernel functions. The kernels proposed in Section \ref{sec:kernels} can play an important role here. Alternatively, one can consider the kernels based on structured features, such as the Gaussian kernel $K(\mathbf{x_i},\mathbf{x_j})=\exp(-\gamma |{\mathbf{x_i-x_j}}|^2),\gamma>0$, where the feature space is implicit and high dimensional. The classifier formed is the radial basis function (\textbf{RBF}) \textbf{SVM}, where the separating hyperplane involving $f(\mathbf{x})=\sum_{i=1}^{N}\alpha_i y_i \,\exp(-\gamma |{\mathbf{x-x_i}}|^2)+b$. Here, $\gamma$ is used as a tuning parameter to determine the smoothness of the decision boundary formed. A large value of $\gamma$ used corresponds to a more complex boundary formed that is more prone to overfitting and vice versa. Depending on the complexity and the nature of the classification problem, different kernels or classifier types may be found to be more suitable.

There are various available software packages that support various types of (regularised) SVM. In our study, the \texttt{R} packages used are \texttt{LiblineaR} and \texttt{e1071}.

The \texttt{LiblineaR} package is particularly useful and efficient for large data sets (especially with many features).
In subsequent experiments, only the following three types of support vector classifiers (SVCs) using primal solvers will be used as the equivalence between primal and dual solvers was shown in \citep{fan2008}.
\begin{enumerate}
	\item The L2-regularised logistic regression (\textbf{LiblineaR-s0}) that solves
	$$\underset{\mathbf{w}}\min\quad  \frac{1}{2}\mathbf{w}^T\mathbf{w} +C\sum\limits_{i=1}^l\log(1+e^{-y_i\mathbf{w}^T\mathbf{x_i}});$$
	\item The L2-regularised L2-loss SVC (\textbf{LiblineaR-s2}) that solves
	$$\underset{\mathbf{w}} \min\quad\frac{1}{2}\mathbf{w}^T\mathbf{w} +C\sum\limits_{i=1}^{l}(\max(0,1-y_i\mathbf{w}^T\mathbf{x_i}))^2;$$
	\item The L1-regularised L2-loss SVC (\textbf{LiblineaR-s5}) that solves
	$$\underset{\mathbf{w}}\min\quad  {\mathbf{\left | w \right |}}_{1} +C\sum\limits_{i=1}^l(\max(0,1-y_i\mathbf{w}^T\mathbf{x_i}))^2,$$
	where ${\left|\cdot\right|}_{1}$ is the L1-norm. This type of regularised classifier creates sparse solutions with selection of nonzero features.
\end{enumerate}
For these solvers, the one-vs.-rest approach will be used for multi-class classification. The tuning parameter $C$ in the formulations above helps to weigh between the regularisation of the weights used (solution for $\mathbf{w}$) and the loss function. The former is favoured when $C$ is small or when the number of observations is much less than the number of features (high-dimensional data).

\subsubsection{Tree-based methods}\label{Tree}
A decision tree is a process of searching through then dividing the predictor space (set of possible values for each feature $x_1,x_2,..., x_p$ into distinct, non-overlapping regions). The structure of a classification tree is that splits are made at one of the $p$ predictors and terminates at nodes with the output class predictions for the region. The same prediction will be made for every observation that lies in the same region. The tree is grown by recursive binary splitting with the aim of minimising resultant misclassification error, Gini index, or cross-entropy.

This simple algorithm, that splits the original predictor space into rectangles, is widely known to have easy interpretation but suffers from low accuracy when it begins to overfit the data. The size of the tree is the main tuning parameter, where a large tree performs well in the training set but risks ``overfitting" the data (resulting in high variance, low bias) while a smaller tree involves lower variance but higher bias. Possible strategies include growing the tree to a large size before pruning it into a smaller tree with cross validation (\textbf{Prunned tree}, see \citep{BFOS84}) or using bagging approach to reduce variance from multiple (uncorrelated) decision trees (\textbf{Random forest}, see \citep{Breiman:2001}) or using boosting to sequentially improve the previous models by adding weak learners to correct classification errors (\textbf{Boosted trees}, see \citep{Freund1997,GBDT:2001}).

Key advantages of tree-based methods is that data pre-processing is not required as inputs are directly used to generate the decision trees. Large amounts of input features can be used without modification or scaling with an inherent ``feature selection" process by selecting order of splits according to variable importance. The most important predictors are first chosen for the initial splits in the tree while unimportant features can be left unused. This is particularly advantageous in our context as it is able to help us interpret the uses of the binned features and identify the events (birth, death or persistence) that are critical for good predictions. The variable importance measures allow us to see which are the important predictors.

\subsubsection{Neural Networks with Dropout}\label{Neural Network}
Neural Networks (NNs) are a popular choice of mathematical modelling inspired by biological neural networks. It works on the basis of transforming the linear combinations of original input features through the use of activation functions and has the effect of creating highly non-linear classifiers. As such, the method is known for its high accuracy and predictive abilities, which often surpasses robust methods like SVM. The complexity of the model depends on the number of hidden layers and the number of units in each layer. Information is passed down from the input layer to the hidden layers to the final output layer through the use of activation functions and linear combinations of weights on units. The output function for binary classification is the sigmoid function while that for multi-class classification is the softmax function.

Many extensions of NNs such as convolution neural networks (CNNs) have been employed for topological data, including the recent paper \citep{cang:2017topologynet}, which uses multi-task multichannel topological convolutional neural networks to learn from small and noisy training data. Deep NNs involve the embedding of multiple non-linear hidden layers that allows them to model complex relationships between features and output. However, the highly non-linear classifier trained may suffer from the effects of overfitting the training data, resulting in poor out-of-sample test performance. As such, regularisation methods, especially dropout \citep{srivastava2014dropout}, have been introduced to address the overfitting problem. This concept has been widely applied and cited to have obtained surprisingly excellent results with significant improvement to existing standard neural networks.

The key features in the use of dropout is that units are randomly dropped along with their incoming and outgoing connections from the network. Each unit is dropped out randomly and independently of other units. Let $p$ to be the probability of a unit being dropped. The dropped out unit is then temporarily excluded from the network and their corresponding weights will not be updated in that iteration. The method of obtaining weights is still unchanged in a basic feed-forward neural networks setting (through backward propagation and convergence to obtain best weights with lowest cost). The resulting thinned neural network can then achieve the effect of preventing complex co-adaptation in the training of the weights compared to a fully connected neural network by adding noise to the hidden units. At test time, the network is all fully connected and the trained weights are then multiplied with a probability of $1-p$. For more details about the dropout, we refer the readers to \citep{srivastava2014dropout}.

\section{Protein Secondary Structure Classification using PHML}\label{sec:numericalstudy}

\subsection{The Classification Problem}\label{Demo problem}

Our classification problem uses the dataset downloaded from the Structural Classification of Proteins-extended (SCOPe) database \citep{Fox:2014}, from which we randomly sample 450 protein samples for our classification task. 
More details into how the database is structured and various classification tasks can be found in \citep{murzin:1995}.

Each protein sample is labeled as one of three classes, namely 1) all-alpha helix, 2) all-beta pleated sheets, and 3) mixed alpha-beta domains, where the alpha ($\alpha$)-helix and the beta ($\beta$)-pleated sheets are two main types of protein secondary structure and their properties are documented in Appendix \ref{app:ab}. For illustrative purpose, we consider a balanced class problem with an equal number of 150 protein samples from each of the three class. Out of 150 samples in each class, 120 of them are used for training while the remaining 30 are used for testing. Our objective is to classify each test sample into one of these three classes based on its structure information. We consider the coarse-grained (CG) representation, i.e., one alpha-carbon ($C_{\alpha}$) for each residue \citep{KLXia:2014c}. As in the first two steps of the general pipeline for PHML, protein structures will be represented by either the Rips complex (RC) or the Alpha complex (AC). 


\subsubsection{Procedure and Evaluation Methodology}
We conduct our studies with the following steps:
\begin{description}
	\item[Step 1] Obtain RC or AC barcodes (Dim-0, 1, 2) for each of the protein samples (see Section \ref{sec:simplicial_complex}).
	\item[Step 2] Generate three sets (Dim-0, 1, 2) of barcodes with the obtained complexes (see Section \ref{sec:PHA}).
	\item[Step 3] Extract topological feature representation (TFR) from each set of barcodes for each protein sample to obtain an overall data matrix with $n$ observations and $p$ features (see Section \ref{sec:feature}). In this study, we consider two TFRs: 1) barcode statistics (BS) and 2) binned features (BF).
	\item[Step 4] Train the model using the data matrix of the training dataset (see Section \ref{sec:MLalgo}) and evaluate the trained classifier with the test dataset. 
\end{description}

In Step 3, for the TFR of BF, the bin size used was fixed at 0.5\AA~ for both the RC and AC filtrations. Moreover, the maximum filtration distance $L$ was set to be $L=20$\AA~ for RC filtration and $L=50$\AA~ for AC filtration. Notice that there are three dimensions for the protein's structure. Hence, the training data matrix has $n=360$ rows/observations and the number of its columns/features $p$ is equal to 13 (for BS with RC or AC), 360 (for BF with RC), or 900 (for BF with AC).

In Step 4, the common metric for evaluating the performance of models is the overall accuracy of the predictions for the 90 held-out test samples.
However, one can also cross-tabulate the output class predictions with the true classes to evaluate class-specific prediction accuracy. While the use of accuracy measure is commonly used to evaluate the performance of balanced-class classifiers, it is unable to distinguish between the severity of the errors. For example, in the context of our classification problem, the failure to distinguish between class 2 (all-alpha) and class 3 (all-beta) should be considered more severe than errors involving class 1 (mixed alpha and beta). Therefore, one can further distinguish the types of misclassifications into 2 categories:
\begin{description}
	\item[Type-I Error] Misclassification errors rate involving class 1 (between classes 1 and 2 and between classes 1 and 3);
	\item[Type-II Error] Misclassification errors rate between classes 2 and 3.
\end{description}
We expect our classifiers to have zero or extremely low second-type error and an acceptable first-type error. For each machine learning model, the output tables generally admit the following format:
\begin{itemize}
	\item 6 Rows: The evaluation criteria to be compared: three individual accuracies of the classes, two types of errors defined above, and the overall accuracy, respectively.
	\item 4 Columns: Corresponding results for the combinations of the types of complex used (RC or AC) and the TFRs used (BS or BF).
\end{itemize}


\subsection{Results} \label{sec:results}

\subsubsection{Results using Support Vector Machine}\label{Results using Support Vector Machine}

Besides of the three types of SVCs in Section \ref{svm}, we also add the RBF SVM into the comparison. We rely on \texttt{R} packages \texttt{e1071} and \texttt{LiblineaR} to implement these four classifiers. Hyperparameters were searched from a grid of $C: 2^{-14},2^{-13},\cdots,2^{14}$ while $g: 2^{-6},2^{3}$ (only for kernels). Cross-validation was performed using the built-in functions in the packages to select optimal parameters ($C$ and ${\gamma}$ for RBF SVM classifiers, ${C}$ for linear SVCs). The above procedure was referenced from that recommended in \citep{Chang:2011}.
\begin{table}[!ht]
	\centering
	\small
	\caption{Results using SVM with different simplicial complexes and TFR methods.}
	\begin{tabular}{|c||c|c|c|c||c|c|c|c|}
		\hline
		\multirow{2}{*}{\textbf{\makecell{Accuracy/\\Results}}} &  \multicolumn{4}{c||}{\textbf{RBF SVM}}& \multicolumn{4}{c|}{\textbf{LiblineaR-s0}} \\
		\cline{2-9}
		&RC-BS & RC-BF & AC-BS & AC-BF
		&RC-BS & RC-BF & AC-BS & AC-BF \\ \hline
		\makecell{Mixed\\Alpha-Beta} &
		83.3 \% & 73.3\% & 80.0 \% & 73.3\% &
		70.0\% & 66.7\% & 43.3\% & 60.0\% \\ \hline
		\makecell{All\\Alpha} &
		86.7 \% & 93.3\% & 86.7 \% & 93.3\% &
		73.3\% & 93.3\% & 93.3\% & 96.7\%  \\ \hline
		\makecell{All\\Beta} &
		86.7 \% & 73.3\% & 90.7 \% & 83.3\% &
		86.7\% & 100\% & 96.7\% & 100\%  \\ \hline
		\makecell{Type-I\\Error} &
		13/90 & 18/90 & 13/90 & 15/90  &
		18/90 & 12/90 & 19/90 & 13/90   \\ \hline
		\makecell{Type-II\\Error} &
		0/90 & 0/90 & 0/90 & 0/90  &
		0/90 & 0/90 & 1/90 & 0/90   \\ \hline
		\makecell{Overall} &
		85.6 \% & 80.0\% & 85.6 \% & 83.3\% &
		80.0\% &  86.7\% & 77.8\% & 85.6\%  \\ \hline \hline
		\multirow{2}{*}{\textbf{\makecell{Accuracy/\\Results}}}&
		\multicolumn{4}{c||}{\textbf{LiblineaR-s2}}&
		\multicolumn{4}{c|}{\textbf{LiblineaR-s5}}\\
		\cline{2-9}
		&RC-BS & RC-BF & AC-BS & AC-BF
		&RC-BS & RC-BF & AC-BS & AC-BF \\ \hline
		\makecell{Mixed\\Alpha-Beta} &
		66.7\% & 66.7\% &23.3\% & 56.7\% &
		66.7\% & 56.7\%  & 26.7\% & 63.3\% \\ \hline
		\makecell{All\\Alpha} &
		83.3\% & 93.3\% &96.7\% & 96.7\% &
		83.3\% & 93.3\% & 26.7\% & 96.7\% \\ \hline
		\makecell{All\\Beta} &
		86.7\% &96.7\% & 100\% & 96.7\% &
		86.7\% & 76.7\% & 96.7\% & 100\% \\ \hline
		\makecell{Type-I\\Error} &
		19/90 & 12/90 & 24/90 & 14/90  &
		19/90 & 22/90 & 23/90 & 12/90  \\ \hline
		\makecell{Type-II\\Error} &
		0/90 & 1/90 & 1/90 & 0/90  &
		0/90 & 0/90 & 1/90 & 0/90 \\ \hline
		\makecell{Overall} &
		78.9\% & 85.6\% & 72.2\% & 84.4\% &
		78.9\% & 75.5\% & 73.3\% & 86.7\% \\  \hline
	\end{tabular}

	\label{tab:Task1_SVM}
\end{table}

Table \ref{tab:Task1_SVM} presents the results for the protein structure classification using SVM. From the results, we can observe that if we adopt BS for TFR, \textbf{RBF SVM} uniformly outperform other SVCs, implying that the data represented by BS are not linearly separable and transformation with Gaussian kernel is helpful; if we adopt BF for TFR, the SVCs generally outperform \textbf{RBF SVM} and among them, L2-regularised logistic regression (\textbf{LiblineaR-s0}) and L2-regularised L2-loss SVC (\textbf{LiblineaR-s2}) perform similarly while L1-regularised L2-loss SVC (\textbf{LiblineaR-s5}) performs the best when we consider AC and BF. The Type-II error occurs only when we consider \textbf{LiblineaR-s2} with RC and BF and when AC and BS are used, whom we do not recommend. We postulate that among all combinations of the types of complex used and the TFRs used, AC and BF carry the most structure information; however, they brought the challenge of high-dimensional statistics as $n\approx p$. SVCs with regularizations can remedy the overfitting issue and give a decent classification result.

\subsubsection{Results using Tree-based methods}\label{Results using Tree-based methods}
A baseline of comparison across tree-based methods was set using a single \textbf{Pruned tree} implemented using the \texttt{R} package \texttt{tree}. A 10-fold (default value) cross-validation was performed to search for the optimal size of tree. \textbf{Random forest} was implemented with $mtry=\sqrt{p}$ (number of features randomly sampled as candidates at each split) and $n_{trees}=1000$ (number of bootstrapped trees to grow) using the \texttt{R} package \texttt{randomForest}. \textbf{Boosted trees} were implemented using \texttt{R} package \texttt{adaBag} with a detailed accompanying tutorial \citep{Alfaro2013}. The boosted trees used are the stumps (with maxdepth=1) and a fixed maximum number of iterations/trees (e.g. 100, 200) were used before searching for the optimal number of iterations that prevents overfitting of the errors. The optimal number of iterations is found by taking value that corresponds to lowest test error. Signs of overfitting are observed when the test error starts to increase after a certain number of iteration and the enssemble can be pruned as well using the built-in function. The package is also suitable for multi-class classification problems where AdaBoost-SAMME \citep{EMSL:2009} is adopted.

\begin{table}[!ht]
	\centering
	\small
\caption{Results using tree-based methods with different simplicial complexes and TFR methods. *Tuned parameter refers to the size of pruned tree or the number of iterations in the boosting approach.}
	\begin{tabular}{|c||c|c|c|c||c|c|c|c|}
		\hline
		\multirow{2}{*}{\textbf{\makecell{Accuracy/\\Results}}} &  \multicolumn{4}{c||}{\textbf{Pruned tree}}&
		\multicolumn{4}{c|}{\textbf{Boosted trees}} \\
		\cline{2-9}
		&RC-BS & RC-BF & AC-BS & AC-BF
		&RC-BS & RC-BF & AC-BS & AC-BF \\ \hline
		\makecell{Mixed\\Alpha-Beta} &
		70.0 \% & 63.3\% & 36.7 \% & 76.7\% &
		76.7\% & 80.0\% & 66.7\% & 96.7\% \\ \hline
		\makecell{All\\Alpha} &
		90.0 \% & 93.3\% & 86.7 \% & 93.3\% &
		86.7\% & 93.3\% & 83.3\% & 86.7\% \\ \hline
		\makecell{All\\Beta} &
		90.0 \% & 86.7\% & 80.0 \% & 76.7\% &
		93.3\% & 93.3\% & 83.3\% & 96.7\% \\ \hline
		\makecell{Type-I\\Error} &
		15/90 & 17/90 & 26/90 & 16/90    &
		13/90 & 10/90 & 20/90 & 6/90  \\ \hline
		\makecell{Type-II\\Error} &
		0/90 & 0/90 & 3/90 & 0/90  &
		0/90 & 0/90 & 0/90 & 0/90 \\ \hline
		\makecell{Overall} &
		83.3 \% & 81.1\% & 67.8 \% & 82.2\% &
		85.6\% & 88.9\% & 77.8\% & 93.3\% \\ \hline
		\makecell{Tuned*\\parameter} &
		3 & 9 & 7 & 4 &
		9 & 37 & 74 & 81  \\ \hline \hline
		\multirow{2}{*}{\textbf{\makecell{Accuracy/\\Results}}} & \multicolumn{4}{c||}{\textbf{Random forest}} \\ \cline{2-5}
		&RC-BS & RC-BF & AC-BS & AC-BF \\ \cline{1-5}
		\makecell{Mixed\\Alpha-Beta}  &
		86.7\% & 86.7\% & 66.7\% & 73.3\% \\ \cline{1-5}
		\makecell{All\\Alpha} &
		86.7\% & 93.3\% & 83.3\% & 90.0\% \\ \cline{1-5}
		\makecell{All\\Beta} &
		93.3\% & 86.7\% & 83.3\% & 96.7\% \\ \cline{1-5}
		\makecell{Type-I\\Error} &
		10/90 & 10/90 & 20/90 & 12/90 \\ \cline{1-5}
		\makecell{Type-II\\Error} &
		0/90 & 0/90 & 0/90 & 0/90  \\ \cline{1-5}
		\makecell{Overall} &
		88.9\% & 88.9\% & 77.8\% & 86.7\% \\ \cline{1-5}
	\end{tabular}
	
	\label{tab:Task1_Tree}
\end{table}

Table \ref{tab:Task1_Tree} reports the results using tree-based methods. Except for \textbf{Pruned tree} with AC and BS, there were no misclassifications between class 2 and 3 (Type-II error). Using AC and BS still yield the worst results as with SVM methods. No matter what complex or TFR are used, bagging (\textbf{Random forest}) and boosting (\textbf{Boosted trees}) can significantly improve the out-of-sample performance. Again, the best result can be found at using AC and BF together with \textbf{Boosted trees} and the best accuracy could hit 93.3\%.

With the careful selection of appropriate tuning parameters such as the number of splits/depth and the number of iterations, modified tree-based models can prevent over-fitted errors or generalize well by converting multiple weak learners into a strong learner. Moreover, tree-based methods possess the advantage of dealing with many correlated input features as there is an inherent ``variable selection" process when searching through all variables for the best split. Using variable importance measures, one can also infer the key barcodes to distinction. Aggregating the top 5 variable importance results from all boosting and random forest models, we notice that none of them are binned features involving Dim-0 as expected and the important variables are consistently identified from the same binned barcodes except for its order of importance.

\subsubsection{Results using Neural Networks}\label{Results using Neural Networks}

The results for types of neural networks are reported: 1) the simplest one-hidden-layer NN; 2) a multiple-layer NN with dropout. For the first type of NN, we simply apply \texttt{R} package \texttt{nnet} but we repeated the training 30 times to obtain the best weights with the lowest cost to avoid a solution from converging to a local minimum. We call it \textbf{Repeated Net}. For the second type of NN, the dropout network was implemented using the \texttt{R} package \texttt{keras}.

In our application, the number of training parameters is often much greater than the number of the available samples for the use of (deep) NNs, which possibly result in overfitting. The data pre-processing, involved removal of constant variables with near-zero variance and subsequent normalisation, is required for the implementation. Preliminary tests without such pre-processing resulted in poor training of weights as the variables with greater numeric ranges dominate those in smaller numeric ranges. The resulting model has even poor training accuracy whereby the in-sample error is inconclusive as well.

Notice that there is no formula for the number of hidden layers and the number of units in the hidden layer. One can choose to add or remove layers conveniently using the \texttt{keras} package. Cross-validation can be used to tune for the optimal number of units but it also involves a series of trial-and-error and is heavily dependent on the nature of the data as well. In a recent paper \citep{Lin2017}, it concludes that a neural network using $p$ input features should not be fitted with fewer than $U=2^p$ units in the hidden layer. However, due to limited computational time and power, the number of units in the single hidden layer was set at $U=2^7$ (resp. $2^4$) for all experiments involving \textbf{Repeated Net} and BS (resp. BF) input features. Errors involving insufficient memory allocation might be presented when the number of hidden units is too large.

Deep neural networks (\textbf{Deep Net}) with or without dropout were applied using the \texttt{R} package \texttt{keras}. The network architecture used is ``$p$ (Input)-$X$-$X$-2$X$-3 (Output)," where $X=128$ refers to the number of units in the first two hidden layers while twice that is used in the last hidden layer. The size of each hidden layer used in the \textbf{Deep Net} is modified with reference to that used in \citep{srivastava2014dropout}. The dropout approach with $p=0.3$ was performed for all hidden layers. The network with Rectified Linear Units (ReLUs) was trained using stochastic gradient descent with a learning rate of 0.4, a decay rate of 0.01, an epoch size of 450, a batch size of 200, and validation split of 20\%.

\begin{table}[!ht]
	\centering
	\small
	\caption{Results using NN with different simplicial complexes and TFR methods.}
	\begin{tabular}{|c||c|c|c|c||c|c||c|c|}
		\hline
		\multirow{3}{*}{\textbf{\makecell{Accuracy/\\Results}}} &  \multicolumn{4}{c||}{\multirow{2}{*}{\textbf{Repeated Net}}} &  \multicolumn{4}{c|}{\textbf{Deep Net}} \\ \cline{6-9}
		& \multicolumn{4}{c||}{} & \multicolumn{2}{c||}{\textbf{No dropout}} & \multicolumn{2}{c|}{\textbf{With dropout}} \\
		\cline{2-9}
		&RC-BS & RC-BF & AC-BS & AC-BF
		&RC-BF & AC-BF & RC-BF & AC-BF\\ \hline
		\makecell{Mixed\\Alpha-Beta} &
		70.0\% & 76.7\% & 63.3\% & 60.0\% &
		90.0\% & 73.3\%  &
		90.0\% & 76.7\% \\ \hline
		\makecell{All\\Alpha} &
		90.0 \% & 93.3\% & 83.3 \% & 100\% &
		90.0\% & 90.0\%  &
		90.0\% & 90.0\% \\ \hline
		\makecell{All\\Beta} &
		83.3\% & 83.3\% & 83.3 \% & 93.3\% &
		90.0\% & 100\%  &
		93.3\% & 100\% \\ \hline
		\makecell{Type-I\\Error} &
		16/90 & 14/90 & 19/90 & 14/90 &
		9/90 & 11/90 &
		8/90 &  10/90 \\ \hline
		\makecell{Type-II\\Error} &
		0/90 & 0/90 & 2/90 & 0/90 &
		0/90 & 0/90 &
		0/90 & 0/90  \\ \hline
		\makecell{Overall} &
		81.1 \% & 84.4\% & 76.7 \% & 84.4\% &
		90.0\% & 87.8\% &
		91.1\% & 88.9\%  \\ \hline
	\end{tabular}
	\label{tab:Task1_NN}
\end{table}

Table \ref{tab:Task1_NN} reports the results using NNs. For the \textbf{Deep Net}, we do not report the results for using BS to save space, because as with most experiments with the ML algorithms presented so far, using BF is better than using BS. It can be seen that a \textbf{Deep Net} is better than a \textbf{Repeated Net} and the dropout approach can slightly improve the performance of a \textbf{Deep Net}. Moreover, the \textbf{Deep Net} provides a relatively high and robust performance with different complexes.

To wrap up the discussions above, we consistently observe that the BF is a better TFR than the BS while different complexes used will slightly affect the performance. When BF is used, it is crucial for the proposed ML algorithm to possess the ability to address the overfitting problem, especially when AC is also used. In our application, \textbf{Boosted tress} and \textbf{Deep Net with dropout} perform the best and the overall accuracy can reach more than 90\%.

\subsection{Further Discussion}
On top of the findings in Section \ref{sec:results}, we further advance our understandings on PHML by exploring the potential improvements from different aspects. In this section, we focus on the BF input features and investigate whether the bin size of BF will be a significant factor in our application. We also explored the usefulness of principal component analysis (PCA) on relieving the high-dimensional learning problem incurred by BF. However, the answer is negative and thus, to save space, we report the corresponding results in the Appendix \ref{app:PCA} for readers' reference.

\subsubsection{Effects of Increasing Bin Number for Binned features}
In this section,
three different bins sizes are considered: 0.5\AA, 0.25\AA, and 0.01\AA. This corresponds to the number of bins of $N$=40, 80, and 200, respectively, for the RC barcodes and $N$=100, 200, and 500, respectively, for the AC barcodes. For each bin, three features are extracted and concatenated and this process is repeated for each of the three dimensions of complexes. Thus, for a protein sample with $N$ bins, the number of features is $9N$. However, the number of training samples is always 360 in our application.

Ideally, the use of more bins may help to represent the protein sample better by identifying more distinguishing features of the barcode (similar to the concept of resolution). However, we do note that the use of too many bins may result in redundant features when the information represented in adjacent bins are the same. The comparison study in this section is to examine the effects of increasing the bin number on the performance of PHML.
We remark that in most cases, due to our relatively small sample size, the number of features exceeds the number of samples ($p>n$), which poses high-dimensional learning challenge.

We repeat the experiments in Section \ref{sec:results} with different bin numbers. The results are reported in Tables \ref{tab:Task1_BF_SVM}-\ref{tab:Task1_BF_DEEPNET} which are placed in Appendix \ref{app:Bin} for better layout.

From Table \ref{tab:Task1_BF_SVM}, we can observe that increasing the bin numbers does not provide improvement for the SVM approaches and it even causes the RBF SVM to become worse. From the statistical perspective, the linear SVCs are robust to the number of bins/features and provide stable performance.

It can be seen from Table \ref{tab:Task1_BF_Tree} that increasing the bin number does improve the performance of the tree-based methods. With enough bins, the \textbf{Random forest} and the \textbf{Boosted trees} can consistently provide about 90\% overall accuracy without any Type-II error. We remark that implied from our experimental results, the number of iterations in \textbf{Boosted trees} should be cross-validated.

The results in Table \ref{tab:Task1_BF_NNREP} clearly indicate that the improvement of \textbf{Repeated Net} by increasing the bin number is significant. For the experiments with \textbf{Deep Net}, we also adjust the network architectures to reflect the increasing number of bins. It can be seen from Table \ref{tab:Task1_BF_DEEPNET} that the \textbf{Deep Net} with more bins as inputs, a more complex network architecture, and the use of dropout improves the performance. By comparing Tables \ref{tab:Task1_BF_NNREP} and \ref{tab:Task1_BF_DEEPNET}, the largest improvement of 6.7\% increase in overall accuracy was observed for two starred numbers. However, when $N=200$ for RC or 500 for AC, the \textbf{Deep Net} performed unsatisfactorily, which could be attributed to many redundant bins formed and serious overfitting issue occurred.

%

\section{Conclusion}\label{Conclusion}
PH-generated barcodes are useful to discover underlying topological features using only finite metric data. Our study provides a gentle walk-through of various methods of Topological Feature Representations (TFR) namely the use of Barcode Statistics (BS) or binned features (BF) before applying them to various machine learning techniques such as SVM, tree-based methods, and NN. The purpose of this paper is to provide guidelines for new users of PH softwares and to provide connections between each stage of analysis. This exposes new users to a complete experience of starting with raw data inputs to the final stage of using machine learning techniques on TFRs. For more advanced users, various stages of analysis can be further improved and fine-tuned depending on their needs and available resources.

Although we only illustrate the PHML on point cloud data, the PHML procedure summarized in this paper could also be applied other datasets including networks or image datasets; see Figure \ref{fig:data_SimCom}. An interesting future research is (extremely) imbalanced class problems that are prevalent in some classification problems.
From the computational perspective, the methods can also be scaled up for larger datasets by paralleling some processes.


\acks{This research is partially supported by Nanyang Technological University Startup Grants M4081840 and M4081842, Data Science and Artificial Intelligence Research Centre@NTU M4082115, and Singapore Ministry of Education Academic Research Fund Tier 1 M401110000.}


\appendix
\section{Properties of the Protein Secondary Structure} \label{app:ab}
This section gives a review of some of the key properties of the secondary structure of proteins. The two main types of protein secondary structure are the alpha ($\alpha$)-helix and the beta ($\beta$)-pleated sheets.

The \textbf{$\alpha$-helix} has the following properties:
\begin{enumerate}
	\item Bond length between immediate $C_{\alpha}$ atom is 3.8\AA.
	\begin{itemize}
		\item This corresponds to the length of typical Betti-0 (Dim-0) bars.
	\end{itemize}
	\item Each turn is made up of 3.6 amino acid residues.
	\begin{itemize}
		\item The formation of Betti-1 (Dim-1) bars can be explained using the slicing technique described in \citep{KLXia:2014c}.
		\item The alpha-helix structure is stabalised by the presence of hydrogen bonds between the amine hydrogen N-H and carbonyl C$=$O oxygen.
		\item Each set of 4 $C_{\alpha}$ atoms form a one-dimensional loop which contributes to a Betti-1 (Dim-1) bar.
	\end{itemize}
	\item Absence of Betti-2 (Dim-2) bars where no cavity is formed since there is insufficient "time" such that the loops are filled up as faces before the cavity can be formed for a single alpha helix.
\end{enumerate}

The \textbf{$\beta$-pleated sheets} have the following properties:
\begin{enumerate}
	\item Bond length/distance between immediate $C_{\alpha}$ atom in the same strand is also 3.8\AA.
	\begin{itemize}
		\item This corresponds to the length of typical Betti-0 (Dim-0) bars. The bar terminates once the atoms are connected.
	\end{itemize}
	\item Each $\beta$-pleated sheet is a stretched out polypeptide chain made up of 3 to 10 amino acid residues.
	\item The $\beta$-pleated sheets are extended structures that are stabalised by hydrogen bonds between residues in adjacent chains.
	\item Each strand must be connected to adjacent strands where the shortest distance between $C_{\alpha}$ and the nearest neighbour in adjacent strand is 4.1\AA.
	\item Adjacent chains run parallel or antiparallel to one another.
\end{enumerate}

\section{Principal Component Analysis on Binned Features} \label{app:PCA}
In this appendix, we investigate the effects of principal component analysis (PCA) on PHML for our application in Section \ref{sec:numericalstudy}. We do not involve the tree-based methods with PCA because trees process dimension reduction by their own construction.

\begin{figure}[!ht]
	\begin{center}
		\begin{tabular}{c}
			\includegraphics[width=0.7\textwidth]{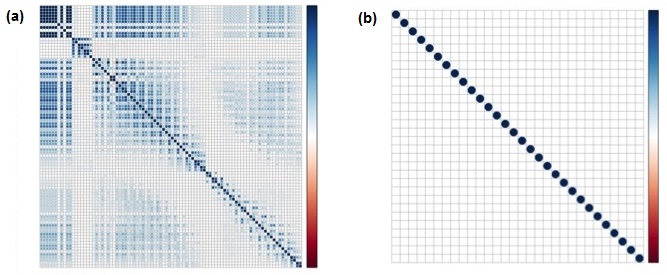}
		\end{tabular}
	\end{center}•
	\caption{({\bf a}) Correlation matrix for the BFs from RC barcodes with 40 bins (where the near-zero variance variables removed). ({\bf b}) Correlation matrix for the 30 PCs transformed from the same BFs. The blue regions corresponds to high positive correlation and is observed between multiple BFs, in contrast with that in the PCs, where white regions correspond to zero correlation.
	}
	\label{fig:Rips20_N40_NZRremoved_Corrplot}
\end{figure}

There are signs of quite high correlation between adjacent BF as seen in Figure \ref{fig:Rips20_N40_NZRremoved_Corrplot}. The use of bins unavoidably suffers from the curse of dimensionality, especially when there are limited number of samples $n$ and $n\ll p$, where $p$ is the number of features. To tackle such a situation, PCA can be applied to transform features into a few uncorrelated PCs, which can be viewed as new features in a lower dimensional feature space. The downside of such an approach is that the final PCs do not have a clear interpretation to the original bins.

In subsequent reports, the experimental results involving principal components transformed from BFs are denoted by an extension ``PC". For consistency, only the first 30 PCs will be used for all transformed features using either RC or AC barcodes. The same set of PCs are used as input features for SVMs and (deep) NNs (with dropout). The settings are the same as specified in Sections \ref{Results using Support Vector Machine} and \ref{Results using Neural Networks}. Tables \ref{tab:Task1_BF_PCA_SVM} and \ref{tab:Task1_BF_PCA_NN} report the corresponding results.

\begin{table}[!ht]
	\centering
\caption{Results using SVM with different simplicial complexes and principal components of BFs.}
	
	\begin{tabular}{|c|c|c|c|c|c|c|c|c|}
		\hline
		\multirow{2}{*}{\textbf{\makecell{Accuracy/\\Results}}}
		&\multicolumn{2}{c|}{\textbf{RBF SVM}} & \multicolumn{2}{c|}{\textbf{LiblineaR-s0}}&
		\multicolumn{2}{c|}{\textbf{LiblineaR-s2}}&
		\multicolumn{2}{c|}{\textbf{LiblineaR-s5}}\\
		\cline{2-9}
		& RC-PC  & AC-PC & RC-PC & AC-PC & RC-PC  & AC-PC & RC-PC & AC-PC\\ \hline
		\makecell{Mixed\\Alpha-Beta} &
		73.3\% & 76.7\%  &
		40.0\% & 36.7\% &
		43.3\% & 36.7\% &
		40.0\% & 40.0\% \\ \hline
		\makecell{All\\Alpha} &
		93.3\% & 93.3\%  &
		96.7\% & 96.7\% &
		96.7\% & 96.7\%  &
		100\% & 96.7\%\\ \hline
		\makecell{All\\Beta} &
		73.3\% & 76.7\%  &
		93.3\% & 96.7\% &
		93.3\% & 96.7\% &
		86.7\% & 96.7\%\\ \hline
		\makecell{Type-I\\Error} &
		18/90 & 16/90   &
		20/90 & 21/90 &
		19/90 & 21/90 &
		22/90 & 20/90 \\ \hline
		\makecell{Type-II\\Error} &
		0/90 & 0/90   &
		1/90 & 0/90 &
		1/90 & 0/90 &
		0/90 & 0/90 \\ \hline
		\makecell{Overall} &
		90.0\% & 73.3\%  &
		90.0\% & 76.7\% &
		90.0\% & 76.7\% &
		75.6\% & 77.8\%\\ \hline
	\end{tabular}
	\label{tab:Task1_BF_PCA_SVM}
\end{table}

\begin{table}[!ht]
	\centering
\caption{Results using NN with different simplicial complexes and principal components of BFs.}
	
	\begin{tabular}{|c|c|c|c|c|}
		\hline
		\multirow{2}{*}{\textbf{\makecell{Accuracy/\\Results}}}   
		&\multicolumn{2}{c|}{\textbf{Repeated NN}} & \multicolumn{2}{c|}{\textbf{\makecell{Deep Networks\\ with dropout}}}\\
		\cline{2-5}
		& RC-PC  & AC-PC & RC-PC & AC-PC \\ \hline
		\makecell{Mixed\\Alpha-Beta} &
		66.7\% & 60.0\%  &
		60.0\% & 56.7\% \\ \hline
		\makecell{All\\Alpha} &
		93.3\% & 96.7\%  &
		93.3\% & 93.3\% \\ \hline
		\makecell{All\\Beta} &
		86.7\% & 86.7\%  &
		93.3\% & 96.7\% \\ \hline
		\makecell{Type-I\\Error} &
		16/90 & 17/90   &
		16/90 & 16/90 \\ \hline
		\makecell{Type-II\\Error} &
		0/90 & 0/90   &
		0/90 & 0/90 \\ \hline
		\makecell{Overall} &
		82.2\% & 81.1\%  &
		82.2\% & 82.2\% \\ \hline
	\end{tabular}
	\label{tab:Task1_BF_PCA_NN}
\end{table}

By comparing the results in Tables \ref{tab:Task1_SVM} and \ref{tab:Task1_BF_PCA_SVM} and Tables \ref{tab:Task1_NN} and \ref{tab:Task1_BF_PCA_NN}, we can see that the use of PCA on BFs does not improve the performance. It implies that the PCA transformation lost information of BFs. In conclusion, it is not recommended that ML algorithms with BFs are incorporated with PCA. However, it does not prohibit PCA from being a powerful visualization tool for the unstructured topological data.

\section{Effects of Increasing Bin Number for Binned features} \label{app:Bin}
In Tables \ref{tab:Task1_BF_SVM}-\ref{tab:Task1_BF_DEEPNET}, the first three or four columns specify the settings of PHML and the remaining columns report the evaluation measurements. The highest overall accuracy number across different bin numbers for a given method is highlighted in red.

\begin{landscape}
\vspace*{\fill}

\begin{table}[htbp]
	\centering
\caption{Results using SVM with different simplicial complexes and BF of different bin numbers.}

	\begin{tabular}{|c|c|c|c|c|c|c|c|c|c|c|c|}
		\hline
		\multirow{2}{*}{\textbf{Steps 1\&2}} & \multirow{2}{*}{\textbf{Step 3}} & \multirow{2}{*}{\textbf{Step 4}}& \multicolumn{3}{c|}{\textbf{$N$=40 (RC) or 100 (AC)}} & \multicolumn{3}{c|}{\textbf{$N$=80 (RC) or 200 (AC)}} &\multicolumn{3}{c|}{\textbf{$N$=200 (RC) or 500 (AC)}} \\
		\cline{4-12}
		& & & Overall & \makecell{Type-I\\Error} & \makecell{Type-II\\Error} & Overall & \makecell{Type-I\\Error} & \makecell{Type-II\\Error}& Overall & \makecell{Type-I\\Error} & \makecell{Type-II\\Error}\\ \hline
		\makecell{Rips \\complex}  & \makecell{BF with \\$L=20$} &RBF SVM &
		\color{red}80.0\% &18/90 &0/90 &
		75.6\% &22/90 &0/90 &
		76.7\% &21/90& 0/90 \\
		\hline
		\makecell{Rips \\complex}  & \makecell{BF with \\$L=20$} &LiblineaR-s0 &
		\color{red}86.7\% &12/90 &0/90 &
		\color{red}86.7\% &12/90 &0/90 &
		82.2\% &15/90& 1/90 \\
		\hline
		\makecell{Rips \\complex}  &\makecell{BF with \\$L=20$} &LiblineaR-s2 &
		\color{red}85.6\%& 12/90 &1/90 &
		83.3\% & 15/90 & 0/90 &
		84.4\% &14/90& 0/90 \\
		\hline
		\makecell{Rips \\complex} &\makecell{BF with \\$L=20$} &LiblineaR-s5 &
		85.6\% &  12/90& 1/90 &
		\color{red}86.7\% &12/90 &0/90 &
		81.1\% &17/90& 0/90\\
		\hline
		\hline
		\makecell{Alpha\\ complex}  & \makecell{BF with \\$L=50$} & RBF SVM &
		\color{red}83.3\% &15/90 &0/90 &
		72.2\% &25/90 &0/90&
		78.9\% &19/90& 0/90 \\
		\hline
		\makecell{Alpha\\ complex} & \makecell{BF with \\$L=50$} &LiblineaR-s0 &
		\color{red}85.6\% &13/90& 0/90&
		84.4\% &14/90& 0/90&
		\color{red}85.6\% &13/90& 0/90 \\
		\hline
		\makecell{Alpha\\ complex}& \makecell{BF with \\$L=50$} &LiblineaR-s2 &
		\color{red}84.4\% &14/90 &0/90 &
		82.2\% &16/90 &0/90&
		82.2\% &16/90& 0/90\\
		\hline
		\makecell{Alpha\\ complex}  & \makecell{BF with \\$L=50$} &LiblineaR-s5 &
		\color{red}86.7\% &12/90 &0/90&
		83.3\% &15/90& 0/90 &
		82.2\% &16/90& 0/90 \\
		\hline
	\end{tabular}
	\label{tab:Task1_BF_SVM}
\end{table}
\vspace*{\fill}
\end{landscape}

\begin{landscape}
\vspace*{\fill}
\begin{table}[htbp]
	\centering
\caption{Results using tree-based methods with different simplicial complexes and BF of different bin numbers. The three maxtree numbers in last two rows correspond to the cross-validated number of iterations for different bin numbers.}

	\begin{tabular}{|c|c|c|c|c|c|c|c|c|c|c|c|}
		\hline
		\multirow{2}{*}{\textbf{Steps 1\&2}} & \multirow{2}{*}{\textbf{Step 3}} & \multirow{2}{*}{\textbf{Step 4}}& \multicolumn{3}{c|}{\textbf{$N$=40 (RC) or 100 (AC)}} & \multicolumn{3}{c|}{\textbf{$N$=80 (RC) or 200 (AC)}} &\multicolumn{3}{c|}{\textbf{$N$=200 (RC) or 500 (AC)}} \\
		\cline{4-12}
		& & & Overall & \makecell{Type-I\\Error} & \makecell{Type-II\\Error} & Overall & \makecell{Type-I\\Error} & \makecell{Type-II\\Error} & Overall & \makecell{Type-I\\Error} & \makecell{Type-II\\Error}  \\ \hline
		\makecell{Rips\\ complex}  & \makecell{BF with\\ $L$=20} & \makecell{Pruned\\ Tree} &
		81.1\% &  17/90& 0/90 &
		83.3\% & 15/90 &0/90&
		\color{red}84.4\% & 14/90 &0/90 \\
		\hline
		\makecell{Alpha\\ complex}  & \makecell{BF with\\ $L$=50} &\makecell{Pruned\\ Tree} &
		82.2\% & 16/90& 0/90 &
		\color{red}84.4\% & 14/90 &0/90&
		80.0\% & 18/90& 0/90 \\
		\hline
		\hline
		\makecell{Rips\\ complex}& \makecell{BF with\\ $L$=20} & \makecell{Random\\ Forest} &
		\color{red}88.9\%  &10/90& 0/90&
		87.8\%  &11/90& 0/90 &
		86.7\%  &12/90 &0/90 \\
		\hline
		\makecell{Alpha\\ complex}  & \makecell{BF with\\ $L$=50} & \makecell{Random\\ Forest} &
		86.7\%  &12/90& 0/90&
		87.8\%  &11/90& 0/90 &
		\color{red}90.0\%  &9/90 &0/90 \\
		\hline
		\hline
		\makecell{Rips\\ complex}  & \makecell{BF with\\ $L$=20} & \makecell{AdaBag \\maxdepth=1\\ maxtree=100} &
		84.4\% & 14/90 &0/90&
		\color{red}88.9\% & 10/90 &0/90&
		\color{red}88.9\%  &10/90 &0/90\\
		\hline
		\makecell{Alpha\\ complex} & \makecell{BF with\\ $L$=50} & \makecell{AdaBag \\ maxdepth=1\\  maxtree=100} &
		84.4\%  &14/90 &0/90&
		85.6\%  &13/90 &0/90&
		\color{red}88.9\%  &10/90 &0/90 \\
		\hline
		\hline
		\makecell{Rips\\ complex}  & \makecell{BF with\\ $L$=20} & \makecell{AdaBag\\ maxdepth=1 \\ maxtree=\\67,54,99} &
		88.9\% & 10/90& 0/90&
		\color{red}90.0\% &9/90 &0/90 &
		\color{red}90.0\% &9/90& 0/90 \\
		\hline
		\makecell{Alpha\\ complex} & \makecell{BF with \\$L$=50} & \makecell{AdaBag\\ maxdepth=1\\maxtree=\\27,11,42}&
		\color{red}93.3\% & 6/90 &0/90 &
		90.0\%  &9/90& 0/90 &
		90.0\% &9/90 &0/90 \\
		\hline
	\end{tabular}
		\label{tab:Task1_BF_Tree}
\end{table}
\vspace*{\fill}
\end{landscape}

\begin{landscape}
\vspace*{\fill}
\begin{table}[htbp]
	\centering
\caption{Results using repeated NN with different simplicial complexes and BF of different bin numbers.}
	\begin{tabular}{|c|c|c|c|c|c|c|c|c|c|}
		\hline
		\multirow{2}{*}{\textbf{Steps 1\&2}} & \multirow{2}{*}{\textbf{Step 3}} & \multirow{2}{*}{\textbf{Step 4}} & \multirow{2}{*}{\makecell{\textbf{Network}\\ \textbf{Architecture}}} &\multicolumn{4}{c|}{\textbf{Accuracy by Class}} & \multicolumn{2}{c|}{\textbf{Error Rates}}\\
		\cline{5-10}
		& & & & Mixed $\alpha$-$\beta$ & All $\alpha$ & All $\beta$ & Overall & Type-I Error & Type-II Error \\ \hline
		\makecell{Rips\\ complex}  & \makecell{BF with\\ $L$=20, $N$=40 }& \makecell{Repeated\\ Net}& 86-16-3&
		76.7\% &  93.3\% &   83.3\%  &84.4\%* & 14/90 & 0/90\\
		\hline
		\makecell{Rips\\ complex} & \makecell{BF with\\ $L$=20, $N$=80} & \makecell{Repeated\\ Net} & 155-16-3&
		70.0\% &  96.7\% &  90.0\% & \color{red}85.6\% & 13/90&  0/90\\
		\hline
		\makecell{Rips\\ complex}  & \makecell{BF with\\ $L$=20, $N$=200} & \makecell{Repeated\\ Net} & 325-16-3&
		73.3\% &  90.0\%&   93.3\%&  \color{red}85.6\% & 13/90  &0/90\\
		\hline
		\hline
		\makecell{Alpha\\ complex}  &\makecell{BF with\\ $L$=50, $N$=100} & \makecell{Repeated\\ Net} & 113-16-3&
		60.0\%&   100\% &  93.3\% & 84.4\% & 14/90&  0/90\\
		\hline
		\makecell{Alpha\\ complex}  &\makecell{BF with\\ $L$=50, $N$=200}& \makecell{Repeated\\ Net} & 207-16-3&
		70.0\%&   93.3\% &  93.3\% & 85.6\% & 13/90 & 0/90\\
		\hline
		\makecell{Alpha\\ complex} & \makecell{BF with\\ $L$=50, $N$=500} & \makecell{Repeated\\ Net} & 513-16-3 &
		80.0\% &  93.3\%  & 93.3\%  &\color{red}88.9\% & 10/90 & 0/90\\
		\hline
	\end{tabular}
	
	\label{tab:Task1_BF_NNREP}
\end{table}
\vspace*{\fill}
\end{landscape}

\begin{landscape}
\vspace*{\fill}

\begin{table}[htbp]
	\centering
\caption{Results using deep NN with different simplicial complexes and BF of different bin numbers.}
	\begin{tabular}{|c|c|c|c|c|c|c|c|c|c|}
		\hline
		\multirow{2}{*}{\textbf{Steps 1\&2}} & \multirow{2}{*}{\textbf{Step 3}} & \multirow{2}{*}{\textbf{Step 4}} & \multirow{2}{*}{\makecell{\textbf{Network}\\ \textbf{Architecture}}} & \multicolumn{3}{c|}{\textbf{No dropout}} & \multicolumn{3}{c|}{\textbf{With dropout}}\\
		\cline{5-10}
		& & & &  \makecell{Overall}& \makecell{Type-I\\Error} & \makecell{Type-II\\Error}& \makecell{Overall} & \makecell{Type-I\\Error} & \makecell{Type-II\\Error} \\ \hline
		\makecell{Rips \\complex}  & \makecell{BF with\\ $L$=20, $N$=40} & Deep Net & 86-128-128-256-3 &
		\color{red}90.0\%  &9/90 & 0/90 &
		\color{red}91.1\%* & 8/90 & 0/90   \\
		\hline
		\makecell{Rips \\complex}  & \makecell{BF with\\ $L$=20, $N$=80} & Deep Net &
		155-256-256-512-3&
		\color{red}90.0\%  &9/90 & 0/90&
		90.0\% & 9/90 & 0/90   \\
		\hline
		\makecell{Rips \\complex}  & \makecell{BF with\\ $L$=20, $N$=200} & Deep Net &
		325-512-512-1024-3&
		81.1\%  &17/90 & 0/90 &
		85.6\% & 13/90 & 0/90  \\
		\hline
		\hline
		\makecell{Alpha \\complex}& \makecell{BF with \\$L$=50, $N$=100} & Deep Net & 113-128-128-256-3 &
		87.8\%   &11/90 & 0/90&
		88.9\%  & 10/90 & 0/90  \\
		\hline
		\makecell{Alpha \\complex}  & \makecell{BF with \\$L$=50, $N$=200} & Deep Net & 207-256-256-512-3 &
		\color{red}88.9\% & 10/90 & 0/90&
		\color{red}91.1\% & 8/90 & 0/90  \\
		\hline
		\makecell{Alpha \\complex}  & \makecell{BF with \\$L$=50, $N$=500 }& Deep Net & 513-512-512-1024-3 &
		86.7\% & 12/90  & 0/90 &
		82.2\% & 16/90 & 0/90  \\
		\hline
	\end{tabular}
	\label{tab:Task1_BF_DEEPNET}
\end{table}
\vspace*{\fill}
\end{landscape}

\vskip 0.2in
\bibliographystyle{plain}
\bibliographystyle{unsrt}

\end{document}